\documentclass[a4paper]{scrartcl}

\usepackage{graphicx,color}
\usepackage{amsmath,amssymb,amsthm}
\usepackage{enumerate}
\usepackage[numbers]{natbib}
\usepackage{rotating}
\usepackage{todonotes}
\usepackage{hyperref}
\usepackage{subcaption}
\usepackage{authblk}

\newtheorem{thm}{Theorem}[section]

\theoremstyle{definition}

\newtheorem{eg}{Example}
\theoremstyle{remark}
\newtheorem{rem}[thm]{Remark}

\newcommand{\tol}{\mathrm{tol}}
\newcommand{\mv}{\mathrm{mv}}
\newcommand{\Pe}{\mathrm{Pe}}
\newcommand{\SR}{\alpha}
\newcommand{\LR}{\nu}
\newcommand{\SI}{\eta}
\newcommand{\LI}{\beta}
\newcommand{\eps}{\varepsilon}

\newcommand{\rme}{\mathrm{e}}
\newcommand{\rmi}{\mathrm{i}}
\newcommand{\CC}{\mathbb{C}}
\newcommand{\RR}{\mathbb{R}}

\newcommand{\argmin}{\mathop{\operatorname{arg\,min}}}
\newcommand{\argmax}{\mathop{\operatorname{arg\,max}}}
\newcommand{\nis}{s} 
\newcommand{\interval}{c} 

\newcommand{\algnt}{\texttt{Alg.1}}

\newcommand{\algre}{\texttt{Alg.2}}
\newcommand{\expmv}{\texttt{expmv}}

\newcommand{\dir}{.}
\graphicspath{{\dir}}

\begin{document}
\title{The Leja method revisited: backward error analysis for the matrix exponential}
\author[1]{M.~Caliari}
\author[2]{P.~Kandolf\thanks{Recipient of a DOC Fellowship of the Austrian Academy of Science at the Department of Mathematics, University of Innsbruck, Austria}}
\author[2]{A.~Ostermann}
\author[2]{S.~Rainer}
\affil[1]{Dipartimento di Informatica, Universit\`{a} di Verona, Italy}
\affil[2]{Institut f\"ur Mathematik, Universit\"at Innsbruck, Austria}
\date{\today}
\maketitle
\begin{abstract}
The Leja method is a polynomial interpolation procedure that can be used to compute matrix functions. In particular, computing the action of the matrix exponential on a given vector is a typical application. This quantity is required, e.g., in exponential integrators.

The Leja method essentially depends on three parameters: the scaling parameter, the location of the interpolation points, and the degree of interpolation. We present here a backward error analysis that allows us to determine these three parameters as a function of the prescribed accuracy. Additional aspects that are required for an efficient and reliable implementation are discussed. Numerical examples illustrating the performance of our Matlab code are included.
\end{abstract}
{\textit{Mathematics Subject Classification} (\textup{2010}): 65F60, 65D05, 65F30}\\
{\textit{Key words:}  Leja interpolation, backward error analysis, action of matrix exponential, exponential integrators, $\varphi$ functions, Taylor series}
\section{Introduction}\label{sec:int}
In many fields of science the computation of the action of the matrix exponential is of great importance. 
As one example amongst others we highlight exponential integrators. 
These methods constitute a competitive tool for the numerical solution of stiff and highly oscillatory problems, see \citep{HO2010}. 
Their efficient implementation heavily relies on the fast computation of the action of certain matrix functions among those the matrix exponential is the most prominent one. 

Given a square matrix $A$ and a vector $v$ the action of the matrix exponential is denoted by $\rme^{{A}}v$.
In general, the exponential of a sparse matrix $A$ is a full matrix. 
Therefore, it is not appropriate to form $\rme^{A}$ and multiply by $v$ for large scale matrices. 
The aim of this paper is to define a backward stable method to compute the action of the matrix exponential based on the Leja interpolation. 
The performed backward error analysis allows one to predict and reduce the cost of the algorithm resulting in a more robust and efficient method. 

For a given matrix $A$ and vector $v$, one chooses a positive integer $\nis$ so that the exponential $\rme^{\nis^{-1}A}v$ can be well approximated. Due to the functional equation of the exponential we can then exploit the relation 
\begin{align}
	\rme^{A}v=\big(\rme^{\nis^{-1}A}\big)^\nis v.
\end{align}
This results in an recursive approximation of $\rme^{A}v=v^{(s)}$ by
\begin{align}\label{eq:scalingsteps}
	v^{(i)}=\rme^{\nis^{-1}A}v^{(i-1)}, \quad v^{(0)}=v.
\end{align}
There are various possibilities to compute the stages $v^{(i)}$. 
Usually, this computation is based on interpolation techniques. 
The best studied methods comprise Krylov subspace methods (see \citep{Sidje1998} and \citep{NiesenWright2012}), truncated Taylor series expansion \citep{Higham2011}, and interpolation at Leja points (see \citep{MR2091132,Leja-comp}).
In this paper we take a closer look on the Leja method (cf.~\eqref{eq:lejamethod} and \eqref{eq:symlejadef} below) for approximating $v^{(i)}\approx~L_{m,\interval}(\nis^{-1}A) v^{(i-1)}$.

Below we present two different ways of performing a backward error analysis of the Leja method. 
Our analysis indicates how the scaling parameter $\nis$, the degree of interpolation $m$ and the interpolation interval $[-\interval,\interval]$ can be selected in order to obtain an appropriately bounded backward error by still keeping the cost of the algorithm at a minimum. 
Furthermore, we discuss how the method benefits from a shift of the matrix and we show how an early termination of the Leja interpolation can help to reduce the cost in an actual computation. 
As a last step we illustrate the stability and behavior of the method by some numerical experiments. 

The paper is structured in the following way. In Section~\ref{sec:be} we introduce the  backward error analysis and draw some conclusions from it. 
In particular we show how this analysis helps us to select the parameters $\nis, m$, and $\interval$. 
In Section~\ref{sec:petc} we discuss some additional aspects for a successful implementation based on the Leja method. 
Section~\ref{sec:numexp} presents some numerical examples dealing with different features and benchmarks for the method. 
In Section~\ref{sec:discussion} we finally give a discussion of the presented results.

For a reader not familiar with the Leja method we included a brief description in Section~\ref{sec:leja}.

\section{The Leja method}\label{sec:leja}

Like every polynomial interpolation, the Leja method essentially depends on the interpolation interval and the position and number of interpolation points. 
The choice of these parameters directly influences the error of the interpolation and the cost. 
In this section we introduce the Leja method based on a sequence of Leja points in a real interval. The extension to a symmetrized complex sequence of points can be found in Section~\ref{sec:complexleja}.

Given an interval $[a,b]$, the Leja points are commonly defined as
\begin{align}\label{eq:rleja}
	\zeta_m\in\argmax_{\zeta\in [a,b]} \prod_{j=0}^{m-1}|\zeta-\zeta_j|, \quad m\geq1,\quad \zeta_0\in\argmax_{\zeta\in [a,b]} |\zeta|.
\end{align}
The interpolation polynomial of the exponential function is then given by 
\begin{align}\label{eq:lejamethod}
	L_m(x;[a,b])=\sum_{j=0}^m \exp[\zeta_0,\ldots, \zeta_j] \prod_{i=0}^{j-1}\left(x-\zeta_i\right),
\end{align}
where $\exp[\zeta_0,\ldots, \zeta_j]$ denotes the $j$th divided difference. 
The scalar interpolation can be extended to the matrix case and rewritten into a two term recurrence relation for the actual computation, see \citep{MR2091132,Leja-comp}.

Due to the functional equation of the exponential it is always possible to shift the argument and perform the interpolation on a symmetric interval with zero as its center. 
This allows us to optimize the algorithm for symmetric intervals.

Let $\zeta_i$ be the Leja points in $[a,b]$ and $\xi_i$ the Leja points in the symmetric interval $[-\interval,\interval]$ with same length.
Then the relation $\zeta_i=\xi_i+\ell$ with $\ell=(a+b)/2$ and $\interval=(b-a)/2$ is valid. 
In practice, we use precomputed points on the interval $[-2,2]$ and scale them to $[-\interval,\interval]$. 
Due to the functional equation of the exponential function the shift can be singled out of the divided differences and we get
\begin{align*}
	L_m(x;[a,b])&=\sum_{j=0}^m \exp[\zeta_0,\ldots, \zeta_j] \prod_{i=0}^{j-1}\left(x-\zeta_i\right) \\&=
	\sum_{j=0}^m \rme^{\ell}\exp[\xi_0,\ldots, \xi_j] \prod_{i=0}^{j-1}\left((x-\ell)-\xi_i\right)\\
	&=\rme^{\ell}L_m(x-\ell;[-\interval,\interval]).
\end{align*}
Therefore, it is always possible to interpolate on a symmetric interval around zero and apply the appropriate shifts to the argument and solution, respectively. 
In the following we will always select $\xi_0=-\interval$ and consequently we get $\xi_1=\interval$ and $\xi_2=0$.
We denote the Leja interpolation polynomial of degree $m$ on the interval $[-\interval,\interval]$ interpolating the exponential by 
\begin{align}\label{eq:symlejadef}
L_{m,\interval}(x) = L_m(x;[-\interval,\interval]).
\end{align} 
Note that it is not necessary to shift the matrix in order to use a symmetric interval. 
Nevertheless, a well chosen shift can lead to faster convergence and can help to avoid round-off and overflow errors.

In order to determine a possible shift we define a rectangle $R=[\SR,\LR]+\mathrm{i}[\SI,\LI]$ in the complex plane.
We do this by splitting up the matrix into its Hermitian part $A_\mathrm{H}$ and skew Hermitian part $A_{\mathrm{SH}}$. 
Furthermore, we find bounds for the field of values and the eigenvalues of these matrices with the help of Gerschgorin's disk theorem, i.e.
\begin{align*}
	\sigma(A)&=\sigma(A_{\mathrm{H}} + A_{\mathrm{SH}})\subseteq \mathcal{W}(A_\mathrm{H} + A_\mathrm{SH})\subseteq \mathcal{W}(A_\mathrm{H}) + \mathcal{W}(A_\mathrm{SH})\\
					&=\mathrm{conv}(\sigma(A_\mathrm{H})) + \mathrm{conv}(\sigma(A_\mathrm{SH})).
\end{align*}
The four real numbers $\SR, \LR, \SI$, and $\LI$ are chosen to satisfy 
\begin{align}\label{eq:recvalues}
	\sigma(A_\mathrm{H})\subseteq[\SR,\LR] \qquad\text{and}\qquad \sigma(A_{\mathrm{SH}})\subseteq\mathrm{i}[\SI,\LI].
\end{align}
We note that in former versions of the Leja method $\LR$ was always assumed nonpositive and $-\SI=\LI$. 
These restrictions are no longer required here.
In this sense, the method is now more general than previous versions.
With the help of the rectangle $R$ the interpolation interval was chosen in \citep{MR2091132,Leja-comp} as the focal interval of the ellipse with smallest capacity circumscribing $R$. 
Here $R$ is used to determine the type of interpolation (real or complex conjugate Leja points) and a possible shift $\mu\in\mathbb{C}$, see Section~\ref{sec:shift}.

We further note that, as stated in \citep{Reichel1990}, the Leja ordering is of great importance for the stability of the method. 
\citeauthor{Reichel1990} suggests the interpolation interval $[-2,2]$. 
The length of the interpolation interval does not influence the numerical accuracy. 
\citeauthor{Reichel1990} suggests $[-2,2]$ only in order to avoid over- and underflow problems which may arise for large 
interpolation intervals and/or with very large values of the interpolation degree.
In this version of the method we deviate from this choice. 
This is admissible since the largest interpolation interval and the largest used degree do not give rise to over- or underflow problems.

\section{Backward error analysis}\label{sec:be}

This section is devoted to the backward error of the action of the matrix exponential when approximated by the Leja method.
We first focus on the interpolation in a real interval, see Section~\ref{sec:complexleja} for the extension to the complex case.

The concept of backward error analysis goes back to Wilkinson, see \citep{Wilkinson1961}. 
The underlying idea is to interpret the result of the interpolation as the exact solution of a perturbed problem $\rme^{A+\Delta A}v$. 
The perturbation $\Delta A$ is the absolute backward error and we aim to satisfy $\|\Delta A\|\leq\tol\,\|A\|$ for a user given tolerance $\tol$.

The here presented backward error analysis exploits a variation of the analysis given in \citep{Higham2011}. 
For this, we define the set
\begin{align*}
	\Omega_{m,\interval}=\{X\in\CC^{n\times n}\colon \rho(\rme^{-X}L_{m,\interval}(X)-I)<1\},
\end{align*}
where $\rho$ denotes the spectral radius and $L_{m,\interval}$ is the Leja interpolation polynomial of degree $m$ on the symmetric interval $[-\interval,\interval]$, see \eqref{eq:lejamethod} and \eqref{eq:symlejadef}. 
Note that $\Omega_{m,\interval}$ is open in $\CC^{n\times n}$ and contains a neighborhood of $0$ for $m\geq2$, since $\xi_2=0$. 
For $X\in\Omega_{m,\interval}$ we define the function
\begin{align}
	h_{m+1,\interval}(X)=\log(\rme^{-X}L_{m,\interval}(X)).
\end{align}
Here $\log$ denotes the principal logarithm \citep[Thm.~1.31]{Higham2008}. 
As $h_{m+1,\interval}(X)$ commutes with $X$ we get $L_{m,\interval}(X)=\rme^{X+h_{m+1,\interval}(X)}$ for $X\in\Omega_{m,\interval}$. 
By introducing a scaling factor $\nis$ such that $\nis^{-1}A\in\Omega_{m,\interval}$ for $A\in\CC^{n\times n}$ we obtain
\begin{align}
	L_{m,\interval}(\nis^{-1}A)^\nis=\rme^{A+\nis h_{m+1,\interval}(\nis^{-1}A)}=:\rme^{A+\Delta A},
\end{align}
where $\Delta A=\nis h_{m+1,\interval}(\nis^{-1}A)$ is the backward error resulting from the approximation of $\rme^A$ by the Leja method $L_{m,\interval}(\nis^{-1}A)^\nis$.

On the set $\Omega_{m,\interval}$ the function $h_{m+1,\interval}$ has a series expansion of the form
\begin{align}\label{eq:seriesextension}
	h_{m+1,\interval}(X)=\sum_{k=0}^\infty a_{k,\interval}X^k.
\end{align}
In order to bound the backward error by a specified tolerance $\tol$ we need to ensure
\begin{align}\label{eq:backwarderror}
	\frac{\|\Delta A\|}{\|A\|}=\frac{\|h_{m+1,\interval}(\nis^{-1}A)\|}{\|\nis^{-1}A\|}\leq\tol
\end{align}
for a given matrix norm. 
As a consequence of this bound, one can select the scaling factor $\nis$ (always a positive integer) such that \eqref{eq:backwarderror} is satisfied for a chosen degree of interpolation $m$. 

In contrast to \citep{Higham2011} we have the endpoint $\interval$ of the interpolation interval as an additional parameter to $m$ and $\nis$ for our analysis. 
In the following, we are going to introduce two different ways of bounding the backward error.
The first one is closely related to the analysis presented in \citep{Higham2011}. 
We study how \eqref{eq:backwarderror} can be used to select the interpolation parameters when we perform a power-series expansion of the backward error. 
In the second approach we consider a contour integral formulation of $h_{m+1,\interval}$ and estimate the error along the contour. 

\subsection{Power-series expansion of the backward error}\label{sec:cd}
In this section we investigate bounds on the backward error represented by $h_{m+1,\interval}$. 
The analysis is based on a power-series expansion of $h_{m+1,\interval}$.

Starting from \eqref{eq:seriesextension} we can bound $h_{m+1,\interval}(X)$ by
\begin{align}\label{eq:htilde}
	\begin{aligned}
\|h_{m+1,\interval}(X)\|&=\left\|\sum_{k=0}^\infty a_{k,\interval} X^k\right\| 
								\leq \sum_{k=0}^\infty |a_{k,\interval}| \left\|X\right\|^k=:\widetilde{h}_{m+1,\interval}(\|X\|).
	\end{aligned}
\end{align}
By inserting this estimate into \eqref{eq:backwarderror} we get 
\begin{align}\label{eq:backwarerror2}
	\frac{\|\Delta A\|}{\|A\|}=\frac{\|h_{m+1,\interval}(\nis^{-1}A)\|}{\|\nis^{-1}A\|}\leq\frac{\widetilde{h}_{m+1,\interval}(\nis^{-1}\|A\|)}{\nis^{-1}\|A\|}.
\end{align}
Since zero is among the interpolation points for $m\ge2$ we get 
\begin{align*}
  \frac{\widetilde{h}_{m+1,\interval}(\theta)}{\theta}=\sum_{k=1}^\infty |a_{k,\interval}| \theta^{k-1}.
\end{align*}
This is a monotonically increasing function for $\theta\geq 0$.  
Furthermore, for $\interval=0$, the Leja interpolation reduces to the truncated Taylor series at zero. 
We thus have $a_{1,0}=0$ for $m\geq1$. 
The equation 
\begin{align}\label{eq:thetamc_equation}
  \frac{\widetilde{h}_{m+1,\interval}({\theta})}{{\theta}}=\tol
\end{align}
therefore has a unique positive root for $\interval$ sufficiently small.
Henceforth, we will call this root $\theta_{m,\interval}$.
The number $\theta_{m,\interval}$ can be interpreted in the following way: 
for the interpolation of degree $m$ in $[-\interval,\interval]$ the backward error fulfills $\|\Delta A\|\leq\tol \|A\|$, if the positive integer $\nis$ fulfills $\nis^{-1}\|A\|\leq\theta_{m,\interval}$. 
In other words, if the norm of a matrix is smaller than $\theta_{m,\interval}$, the interpolation of degree $m$ with points in $[-\interval,\interval]$ has an error less than or equal to $\tol$. 

In the analysis up to now, we only used that zero is among the interpolation points. 
However, the following discussion requires the sequence of Leja points.

In order to compute $\widetilde{h}_{m+1,\interval}$ in a stable manner we expand $h_{m+1,\interval}$ in the Newton basis for Leja points in $[-\interval, \interval]$ as
\begin{align}\label{eq:leja:seriesextension}
	h_{m+1,\interval}(X)=\sum_{k=m+1}^\infty h_{m+1,\interval}[\xi_0,\ldots,\xi_k]\prod_{j=0}^{k-1}(X-\xi_j I),
\end{align}
where $h_{m+1,\interval}[\xi_0,\ldots,\xi_k]$ denotes the $k$th divided difference. The above series starts with $k=m+1$ as $h_{m+1,\interval}(\xi_j)=0$ for $j=0,\ldots, m$. 
Rewriting this series in the monomial basis we obtain \eqref{eq:htilde} with the according coefficients $a_{k,\interval}$. 
In order to get reliable results for these coefficients we use $300$ digits in the actual computation.

\begin{figure}[ht]\centering
	\includegraphics{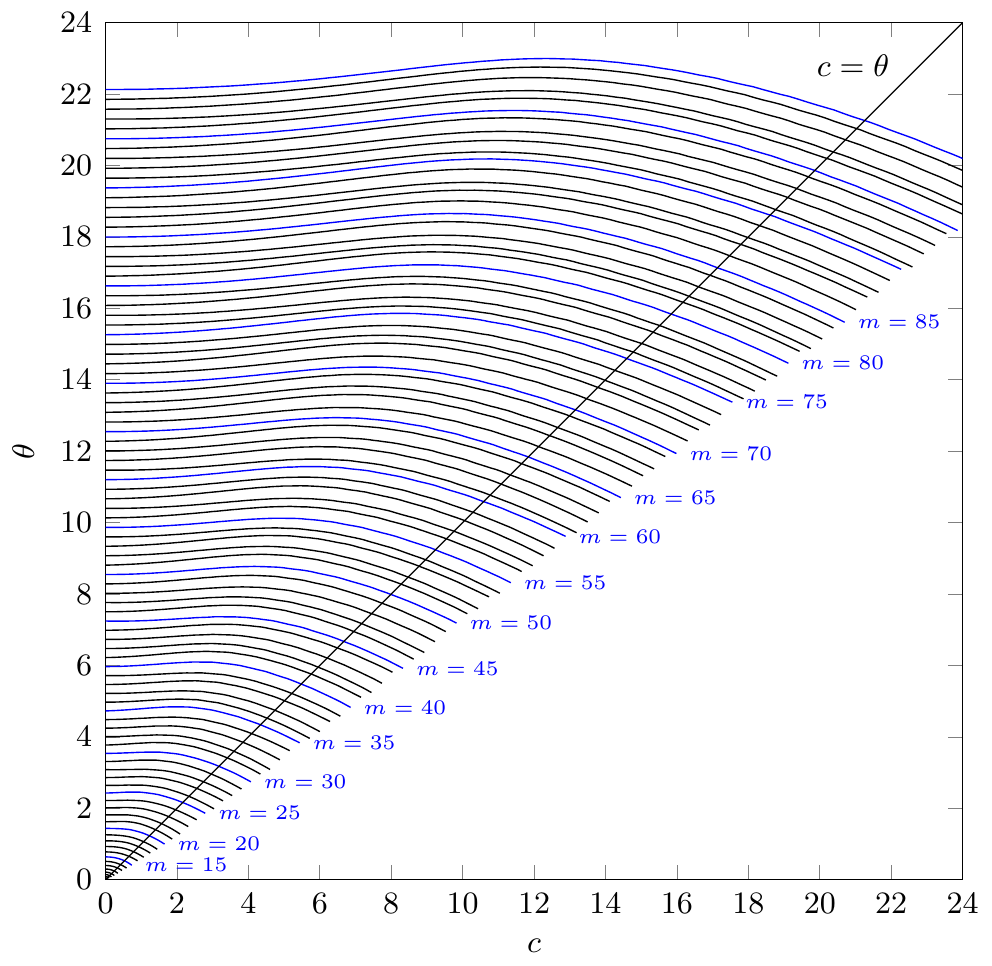}
	\caption{\label{fig:thetavsbeta}The root $\theta=\theta_{m,\interval}$ as a function of the right endpoint $\interval$ of the interpolation interval for real Leja points in $[-\interval,\interval]$. The tolerance is set to $\tol=2^{-53}$. Along each line the interpolation degree $m$ is kept fixed.}
\end{figure}

Figure~\ref{fig:thetavsbeta} displays the path of $\theta_{m,\interval}$ for the Leja interpolation with respect to $\interval$ for fixed $m$ up to $100$ and $\tol=2^{-53}$. 
For an actual computation one has to truncate the series \eqref{eq:leja:seriesextension} at some index $M$. 
We always used $M=3m$.

As we now have a way of bounding the backward error we discuss the choice of the number of scaling steps in \eqref{eq:scalingsteps}. 
We propose to select the integer $\nis$ depending on $m$ and $\interval$ in such a way that the cost of the algorithm becomes minimal. 
We have the limitation that $m$ is bounded by $100$ to avoid problems with over- and underflow. 
However, we get several possibilities to select the free parameter $\interval$ describing the interval. 

The value $\theta_{m,0}$ corresponds to the truncated Taylor series expansion as described in \citep{Higham2011}. 
A second possibility is to choose, for a fixed $m$, $\interval$ in such a way that $\theta_{m,\interval}$ is maximal. 
This corresponds to the interpolation interval that admits the largest norm of $A$. 
A third possibility is to select the interpolation interval such that the right endpoint $\interval$ coincides with $\theta_{m,\interval}$. 
These are the points on the diagonal in Figure~\ref{fig:thetavsbeta}.

A priori none of the above choices can be seen to be optimal for an arbitrary matrix. 
The choice $\interval=0$ together with the smallest $m$ such that $\theta_{m,0}\ge \|A\|$ is a good choice for a matrix $A$ with all the eigenvalues clustered in a neighborhood of $0$.
On the other hand, if the convex hull of the eigenvalues of a matrix $A$, with $\|A\|\approx\delta$, is the interval $[-\delta,\delta]$, the choice $\theta_{m,\delta}$ with smallest $m$ such that $\theta_{m,\delta}\ge \|A\|$ is preferable. 

We choose $\theta_{m,\interval}$ according to the third option, which favors normal matrices where all eigenvalues lie in an interval. 
More precisely, we select
\begin{align}\label{eq:thetam}
	\theta_m = \min\{c \colon \theta_{m,\interval}= c\}.
\end{align}
This means that $\theta_m$ is the first intersection point of the graph $(c, \theta_{m,\interval})$ with the diagonal, cf.~Figure~\ref{fig:thetavsbeta}.

The behavior of the curves $\interval \mapsto (\interval,\theta_{m,\interval})$ is not unexpected.
Let us consider the approximation of  $\mathrm{e}^\theta$ by $L_{m,\interval}(\theta)$ for $\theta\ge 0$ and $\nis=1$.
Then $\tilde h_{m+1,\interval}(\theta)/\theta$ is an overestimate of the relative backward error $h_{m+1,\interval}(\theta)/\theta$.
The value $\theta_{m,\interval}$ represents the maximum value for which $L_{m,\interval}(\theta_{m,\interval})$ is an acceptable approximation of $\mathrm{e}^{\theta_{m,\interval}}$.
The value $\theta_{m,0}$ corresponds to interpolation at a set of confluent points at $\interval=0$, i.e.~the truncated Taylor series approximation.
If we slightly increase the interpolation interval $[-\interval,\interval]$, about half of the interpolation points lie in $[0,\interval]$.
Therefore, it is possible to have an acceptable interpolation up to $\theta_{m,\interval}\ge \theta_{m,0}$.
If we continue to increase $\interval$, the mutual distance between interpolation points increases as well. 
When the interval gets too large, the number of interpolation points is too small to achieve the desired accuracy and the value $\theta_{m,\interval}$ starts to decrease. 

Figure~\ref{fig:thetavsbeta} shows that $\theta_{m,\interval}\geq\theta_m$ for all $0\leq\interval\leq\theta_m$.
Therefore, we can safely interpolate with degree $m$ for all intervals $[-\interval,\interval]$ with $0\leq \interval\leq\theta_m$. 
We will use this fact in Section~\ref{sec:hump}. 
 
We compute $\theta_m$ by a combination of two root finding algorithms based on Newton's method.
The inner equation \eqref{eq:thetamc_equation} for computing $\theta_{m,\interval}$ is solved by an exact Newton iteration with an accuracy of $10^{-20}$. 
The outer equation $\theta_{m,\interval}=\interval$ is also solved by Newton's method.
This time, however, the necessary derivative is approximated by numerical differentiation. 
We compute the result up to an accuracy of $10^{-18}$. 
The resulting values are truncated to 16 digits (double precision) and used henceforth as the $\theta_m$ values. 
In Table~\ref{tab:thetam} we listed selected (rounded) values of $\theta_m$ for various $m$ and certain tolerances. 
\begin{table}\small
\begin{tabular}{r|rrrrrrr}
 \hline
    $m$ &        5 &       10 &       15 &       20 &       25 &       30 &       35 \\
   half & 6.43e-01 & 2.12e+00 & 3.55e+00 & 5.00e+00 & 6.37e+00 & 7.51e+00 & 8.91e+00 \\
 single & 9.62e-02 & 8.33e-01 & 1.96e+00 & 3.26e+00 & 4.69e+00 & 5.96e+00 & 7.44e+00 \\
 double & 1.74e-03 & 1.14e-01 & 5.31e-01 & 1.23e+00 & 2.16e+00 & 3.18e+00 & 4.34e+00 \\
 \hline\hline
    $m$ &       40 &       45 &       50 &       55 &       60 &       65 &       70 \\
   half & 1.00e+01 & 1.10e+01 & 1.23e+01 & 1.35e+01 & 1.48e+01 & 1.59e+01 & 1.71e+01 \\
 single & 8.71e+00 & 1.00e+01 & 1.15e+01 & 1.27e+01 & 1.40e+01 & 1.52e+01 & 1.64e+01 \\
 double & 5.48e+00 & 6.67e+00 & 7.99e+00 & 9.24e+00 & 1.06e+01 & 1.18e+01 & 1.32e+01 \\
 \hline\hline
    $m$ &       75 &       80 &       85 &       90 &       95 &      100 \\
   half & 1.84e+01 & 1.94e+01 & 2.07e+01 & 2.20e+01 & 2.30e+01 & 2.42e+01 \\
 single & 1.76e+01 & 1.87e+01 & 1.99e+01 & 2.12e+01 & 2.23e+01 & 2.35e+01 \\
 double & 1.46e+01 & 1.58e+01 & 1.71e+01 & 1.86e+01 & 1.99e+01 & 2.13e+01 \\
 \hline
\end{tabular}
\caption{\label{tab:thetam}Samples of the (rounded) values $\theta_m$ for tolerances \emph{half} ($\tol=2^{-10}$), \emph{single} ($\tol=2^{-24})$ and \emph{double} ($\tol=2^{-53}$) for the real Leja interpolation.}
\end{table}

We next describe the choice of the parameters used in our implementation. For each $m$ the optimal value of the integer $\nis$ is given by
\begin{align}
	\nis=\lceil\|A\|/\theta_m\rceil.
\end{align}
We recall that we have chosen $m_{\mathrm{max}}=100$.
The cost of the interpolation is dominated by the number of matrix-vector products computed during Newton interpolation. Therefore, the cost is at most 
\begin{align}\label{eq:costfunction}
	C_m(A):=sm=m\lceil\|A\|/\theta_m\rceil,
\end{align} 
resulting in the optimal $m_*$ and corresponding $\nis_*$ and $\interval_*$ as 
\begin{align}\label{eq:scalingstep_circle}
	m_*=\argmin_{2\leq m \leq m_{\mathrm{max}}} \left\{m \left\lceil\frac{\|A\|}{\theta_m} \right\rceil \right\},\quad
	\nis_* = \left\lceil\frac{\|A\|}{\theta_{m_*}} \right\rceil, \quad c_*=\theta_{m_*}.
\end{align} 


\begin{rem}[precomputed divided differences]
We now have a fixed discrete set of interpolation intervals, given by $\theta_m$. Therefore, the associated divided differences can be precomputed, once and for all. 
\end{rem}

\begin{rem}[nonnormal matrices]\label{rem:alphap}
The truncated Taylor series method is able to exploit the values $d_p=\|A^p\|^{1/p}$. As shown in \citep[Eq.~(3.6)]{Higham2011}, the backward error satisfies
\begin{align*}
\frac{\|\Delta A\|}{\|A\|} \le \frac{\widetilde{h}_{m+1,0}(\nis^{-1}\alpha_p(A))}{\nis^{-1}\alpha_p(A)},
\end{align*}
where $\widetilde{h}_{m+1,0} (\theta) = \sum_{k=m+1}^\infty |a_{k,0}|\theta^k$ and $\alpha_p(A) = \min(d_p,d_{p+1})$ with arbitrary $p$ subject to $p(p-1)\le m+1$. 
For nonnormal matrices, this can be a sharper bound as $\alpha_p(A)\ll \|A\|$ is possible. 
Note that our method is not able to use this relation right away. 
This is due to the fact that the series representation of  $\widetilde{h}_{m+1,\interval}$ starts at $k=1$ for $c\ne 0$, see \eqref{eq:htilde}. 
As a result, we might use more scaling steps for such problems, see Section~\ref{sec:numexp} for some experiments. Nevertheless, the values $d_p$ can be favorably used also for the Leja method, see Section~\ref{sec:hump}.
\end{rem}

\subsection{Contour integral expansion of the backward error}\label{sec:convergence_ellipses}

The selection of the scaling step and the length of the interpolation interval based on the norm of the matrix does not take into account the distribution of the eigenvalues. 
In this section we investigate bounds of the backward error, based on a contour integral expansion along ellipses that enclose the spectrum of the matrix.
This introduces more flexibility as an ellipse can vary its shape from an interval to a circle. 
By this we can better capture the distribution of the eigenvalues of a matrix $A$ than by simply taking the number $\|A\|$.

Again our aim is to find a bound for $h_{m+1,\interval}$. 
Due to the fact that zero is among the interpolation points for $m\geq2$ it is convenient to write $h_{m+1,\interval}$ as
\begin{align}
 h_{m+1,\interval}(X)=Xg_{m+1,\interval}(X).
\end{align}
For fixed $\eps>0$ we can rewrite \eqref{eq:backwarderror} in the Euclidean norm as
\begin{align}\label{eq:ellipsestimate}
	\begin{aligned}
	\frac{\|\Delta A\|_2}{\|A\|_2}&=\frac{\|h_{m+1,\interval}(\nis^{-1}A)\|_2}{\|\nis^{-1}A\|_2}\\
			&\leq \|g_{m+1,\interval}(\nis^{-1}A)\|_2 \\
			&= \left\|\frac{1}{2\pi\rmi}\int_{\Gamma}g_{m+1,\interval}(z)(zI-\nis^{-1}A)^{-1}\,\mathrm{d}z\right\|_2\\
			&\leq\frac{\mathcal{L}(\Gamma)}{2\pi\eps}\|g_{m+1,\interval}\|_{\Gamma}.
	\end{aligned}
\end{align}
Here $\Gamma=\partial K$ denotes the boundary of the domain $K$ that contains $\Lambda_\eps(\nis^{-1}A)$, the $\eps$-pseudospectrum of $\nis^{-1}A$. 
The $\eps$-pseudospectrum of a matrix $X$ is given by 
\begin{align*}
 \Lambda_\eps(X)= \left\{z\colon \left\|(zI-X)^{-1}\right\|_2\geq \eps^{-1}\right\}.
\end{align*}
The length of $\Gamma$ is denoted by $\mathcal{L}(\Gamma)$ and $\|\cdot\|_{\Gamma}$ is the maximum norm on $\Gamma$. 
For given $m, \interval$, and $K$ the last term in \eqref{eq:ellipsestimate} can be computed in high precision. 
We use $300$ digits and sample the contour in any performed computation.

For the time being, let us fix $m$. 
Furthermore, we assume that $\Gamma$ is an ellipse, with focal interval equal to the interpolation interval $[-\interval,\interval]$, and its convex hull $K$ encloses $\Lambda_\eps(\nis^{-1}A)$. 
As a result of these assumptions, there is not only one ellipse but rather a two parameter family of ellipses $\Gamma_{\gamma,\interval}$. 
The parameters are the right endpoint $\interval$ of the interpolation (focal) interval and the capacity $\gamma$ of the ellipse, that is the half sum of the semi-axes.
In the following we describe how to extract a discrete set of ellipses from the two parameter family of ellipses. 
This discrete set can then be stored and used in the algorithm.

We start by reducing the two parameter family of ellipses $\Gamma_{\gamma,\interval}$ to a one parameter family, only depending on the focal interval $[-\interval,\interval]$.
For fixed $\interval$ we have a family of confocal ellipses that are described by 
\begin{align}\label{eq:ellipseboundary}
	\Gamma_{\gamma,\interval} = \left\{z\in\CC\colon z=\gamma w + \frac{\interval^2}{4\gamma w},\quad |w|=1\right\}.
\end{align}
An ellipse $\Gamma_{\gamma,\interval}$ will be considered valid for interpolation if 
\begin{align}\label{eq:gammafoc}
	\frac{\mathcal{L}(\Gamma_{{\gamma},\interval})}{2\pi\eps}\|g_{m+1,\interval}\|_{\Gamma_{{\gamma},\interval}}\leq\tol
\end{align}
is satisfied. 
For every $\interval$ there exists an ellipse with largest capacity $\gamma=:\gamma_{m,\interval}$ satisfying \eqref{eq:gammafoc} as tolerance equality, if $m$ is sufficiently large.
We single out this ellipse and thereby link the capacity directly with the focal interval and construct a one parameter family of ellipses.
%

We further reduce the number of ellipses by selecting only a discrete set of focal intervals for every $m$.
More precisely, we use the known values $\theta_j$ for $j\geq m$ from \eqref{eq:thetam} and Table~\ref{tab:thetam}, respectively. 
For these values we already have precomputed divided differences at hand and no extra storage is needed.

The overall procedure is as follows.
For each interpolation (focal) interval $[-\theta_j,\theta_j]$, $j\ge m$ we compute the ellipse with largest capacity fulfilling \eqref{eq:gammafoc} with $\eps=\tfrac{1}{50}$ and store its semi-axes. 
We increase $j$ as long as there is a $\gamma=:\gamma_{m,\theta_j}$ satisfying \eqref{eq:gammafoc}. 
Furthermore, we enforce the upper limit $j\leq120$. 
With this selection we allow at most $20$ ellipses for the maximal degree of interpolation $m_\mathrm{max}=100$.

\begin{figure}[hbt]\centering
	\includegraphics{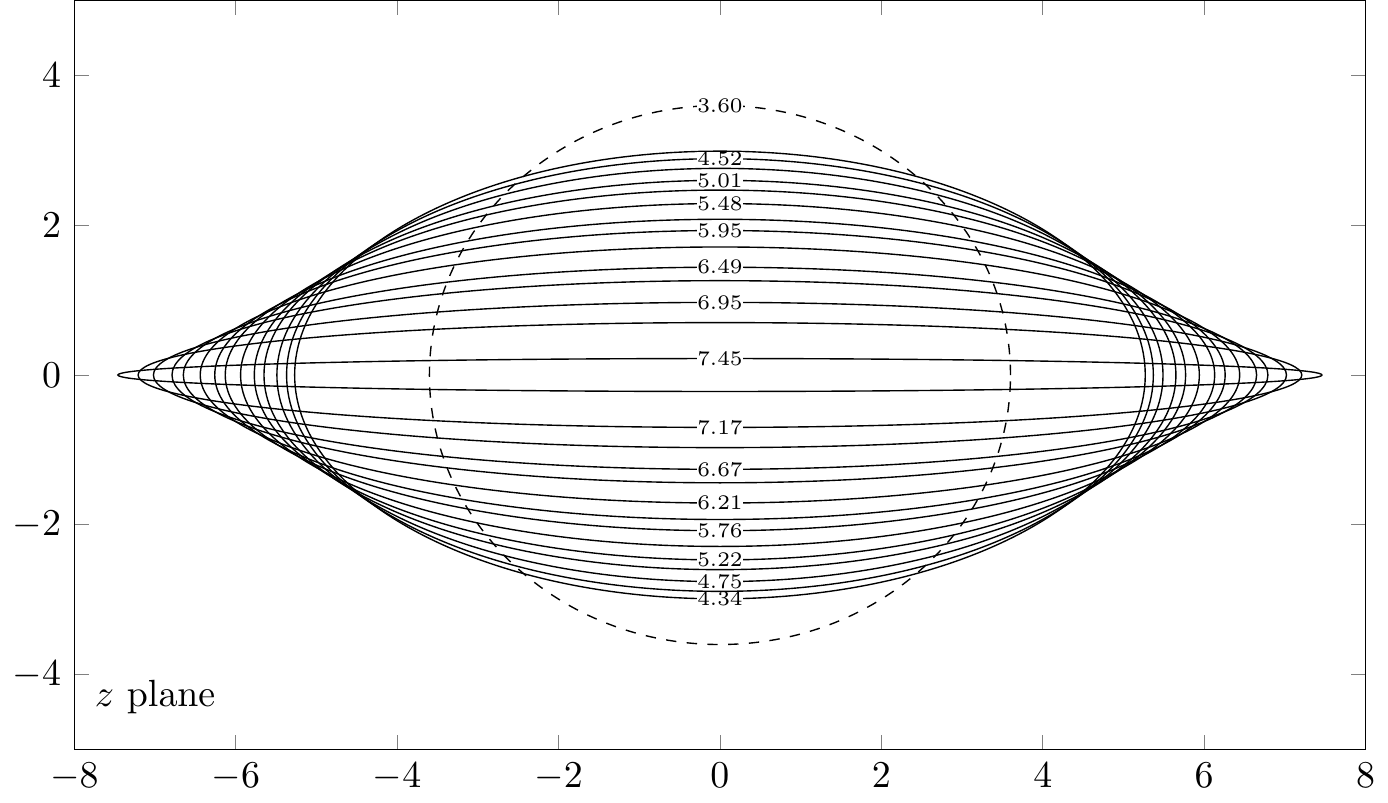}
	\caption{\label{fig:ellips}For $\eps=\tfrac{1}{50}$, $m=35$ and $\tol=2^{-53}$ the ellipses \eqref{eq:ellipseboundary} satisfying \eqref{eq:gammafoc} are shown for various focal intervals $[-\theta_j,\theta_j]$ with $j=35, \ldots, 48$. The value $\theta_j$ (see Table~\ref{tab:thetam}) is indicated on the ellipse. The dashed circle has radius $\theta_{32}=3.60$. It is the largest circle in lieu of the ellipse $\Gamma_{\gamma,\interval}$ that fulfills \eqref{eq:gammafoc} for $m=35$.}
\end{figure}

Figure~\ref{fig:ellips} shows the stored ellipses for $m=35$ and $\tol=2^{-53}$. 
The dashed circle has radius $\theta_{32}=3.60$ corresponding to the largest circle with radius $\theta_j$ that fulfills \eqref{eq:gammafoc} if a circle is used instead of an ellipse; see Section~\ref{sec:hump} for further reasoning why to include this circle. 

We can see that, for larger focal intervals, the semi-minor axis decreases until we end up, in the limit, with an interval on the real axis. 
As we have additional information on the spectrum of the matrix at hand, it is possible to interpolate the exponential of certain matrices with fewer scaling steps than predicted by our power-series estimate.

\begin{rem}If we reduce the size of the focal interval of our ellipses $\Gamma_{m,\interval}$ to a point, we end up with a circle. 
In fact, for a fixed $m$ this circle is slightly smaller than the one obtained by the estimate $\theta_m$. 
\end{rem}

For our backward error analysis we can interpret the value $\gamma_{m,\theta_j}$ in the following way. 
If we prescribe an interval $[-\theta_j,\theta_j]$ and select $m+1$ Leja points in this interval, we have $\|\Delta A\|\leq\tol\,\|A\|$ under the assumption that $\nis\geq1$ is selected such that $\Lambda_\eps(\nis^{-1}A)\subseteq\mathrm{conv}(\Gamma_{\gamma_{m,\theta_j},\theta_j})$. 

\begin{figure}\centering
	\includegraphics{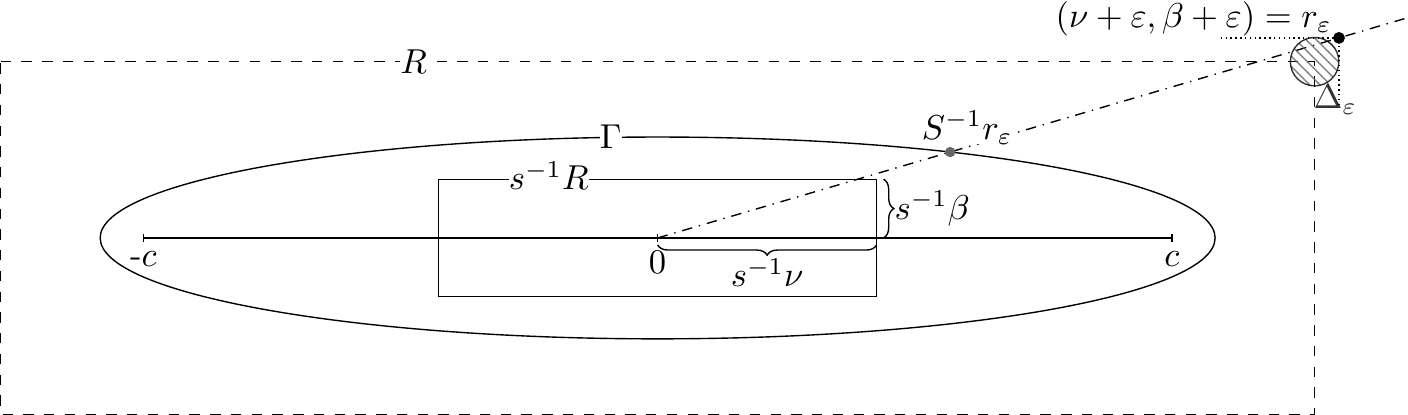}
	\caption{\label{fig:scaling_ellipse}Illustration on the selection of the correct scaling factor $\nis$ to fit the scaled and extended estimate of the pseudospectrum $\nis^{-1}(R+\Delta_\eps)$ inside the ellipse $\Gamma$ with convex hull $K\subset\mathbb{R}^2$. We have $S=\sqrt{(\LR+\eps)^2 a^{-2}+(\LI+\eps)^2 b^{-2}}$ and $\nis=\lceil S\rceil$.} 
\end{figure}

Before discussing how to select the optimal ellipse for $m$, we show how to compute $\nis$ for a given matrix $A$ and an ellipse $\Gamma$. 
We recall that, with the help of Gershgorin's disk theorem, we can enclose the spectrum of $A$ in a rectangle $R$ with vertices $(\SR,\LI), (\SR, \SI), (\LR, \LI), (\LR, \SI)$, see \eqref{eq:recvalues}. 
Furthermore, we assume that this rectangle is centered in zero ($-\SR=\LR$ and $-\SI=\LI$), otherwise we shift the matrix accordingly. 
In order to keep the notation simple we consider a single ellipse $\Gamma$ with focal interval $[-\interval,\interval]$ and capacity $\gamma$ for which \eqref{eq:gammafoc} is satisfied. 
As before we denote the convex hull of $\Gamma$ by $K$.
Let $\Delta_\eps$ denote the open disc of radius $\eps$ around the origin. Hence, we have the following chain of inclusions 
\begin{align*}
	\Lambda_\eps(A) \subseteq \mathcal{W}(A) + \Delta_\eps \subseteq R + \Delta_\eps.
\end{align*}
The first inclusion connecting the pseudospectrum and the field of values can be found in \citep{Trefethen2005}.
The above inclusions state that if $R+\Delta_\eps\subset K$ then $\|\Delta A\|_2\leq\tol\,\|A\|_2$. 
Our aim is to determine the correct scaling factor $\nis$ such that the inclusion $\nis^{-1} (R+\Delta_\eps)\subseteq K$ is valid. 
We do this by computing the intersection of $\Gamma$ with the straight line through zero and $r_\eps=(\LR+\eps,\LI+\eps)$, the upper right vertex of the rectangle extended by $\eps$. The procedure is illustrated in Figure~\ref{fig:scaling_ellipse}. 

For
\begin{align*}
	a=\gamma+ \frac{\interval^2}{4\gamma}, \qquad b=\gamma - \frac{\interval^2}{4\gamma}
\end{align*}
denoting the semi-axes of $\Gamma$ we have 
\begin{align}\label{eq:scaling_single_ellipse}
	\nis=\left\lceil{\sqrt{\frac{(\LR+\eps)^2}{a^2}+\frac{(\LI+\eps)^2}{b^2}}}\right\rceil.
\end{align}
Due to our choice of $r_\eps$ it holds that $\nis^{-1}(R+\Delta_\eps)\subseteq K$ for $\nis^{-1}r_\eps \in K$.

We can now use the the degree of interpolation $m$ to minimize the cost of the interpolation. 
This is done in the following way. 
As discussed above and illustrated in Figure~\ref{fig:ellips}, for every $m$, we get a family of ellipses with semi-axes 
\begin{align}\label{eq:semiaxis}
	a_{m,\theta_j}=\gamma_{m,\theta_j}+ \frac{\theta_j^2}{4\gamma_{m,\theta_j}} 
	\quad \text{and}\quad
	b_{m,\theta_j}=\gamma_{m,\theta_j}- \frac{\theta_j^2}{4\gamma_{m,\theta_j}},
\end{align}
where $\gamma_{m,\theta_j}$ fulfills \eqref{eq:gammafoc}. 
Recall that we have chosen $m_{\mathrm{max}}=100$.
We now use \eqref{eq:scaling_single_ellipse} to select the optimal ellipse for each $m$. 
In this family the optimal ellipse is identified as the one with the fewest scaling steps.
For these optimal ellipses the number of scaling steps $\nis_m$ is given by 
\begin{align*}
	S_{m,j}= \sqrt{ \left(\frac{\LR+\eps}{a_{m,\theta_j}}\right)^2 + \left(\frac{\LI+\eps}{b_{m,\theta_j}}\right)^2}, 
	\quad j_m=\argmin_{j\geq m} \left\{\left\lceil S_{m,j}\right\rceil\right\}, 
	\quad{\nis}_m=\left\lceil S_{m,j_m}\right\rceil.
\end{align*}
Now we can minimize with a similar cost function as in \eqref{eq:costfunction} over $m$ and obtain our optimal degree of interpolation $m_*$ and scaling factor $\nis_*$ as
\begin{align}\label{eq:scalingstep_ellipse} 
	m_*=\argmin_{2\leq m \leq m_{\mathrm{max}}}\left\{m {\nis}_m\right\},\quad \nis_*=\nis_{m_*}.
\end{align} 
The corresponding $c_*$ is given by $\theta_j$ with $j=j_{m_*}$

\subsection{Symmetrized complex Leja points}\label{sec:complexleja}

All the statements made in the previous sections remain valid if we use complex conjugate Leja points \citep{Mashfree}. 
The advantage of such points lies in the better handling of matrices that have eigenvalues with dominant imaginary parts.
This situation is characterized by a height-to-width ratio of more than one for the rectangle $R$. 
Examples include the (real) discretization matrices of transport equations or the discretization of the Schr\"odinger operator (a complex matrix) which has eigenvalues on the negative imaginary axis.
\begin{figure}[th]\centering
	\includegraphics{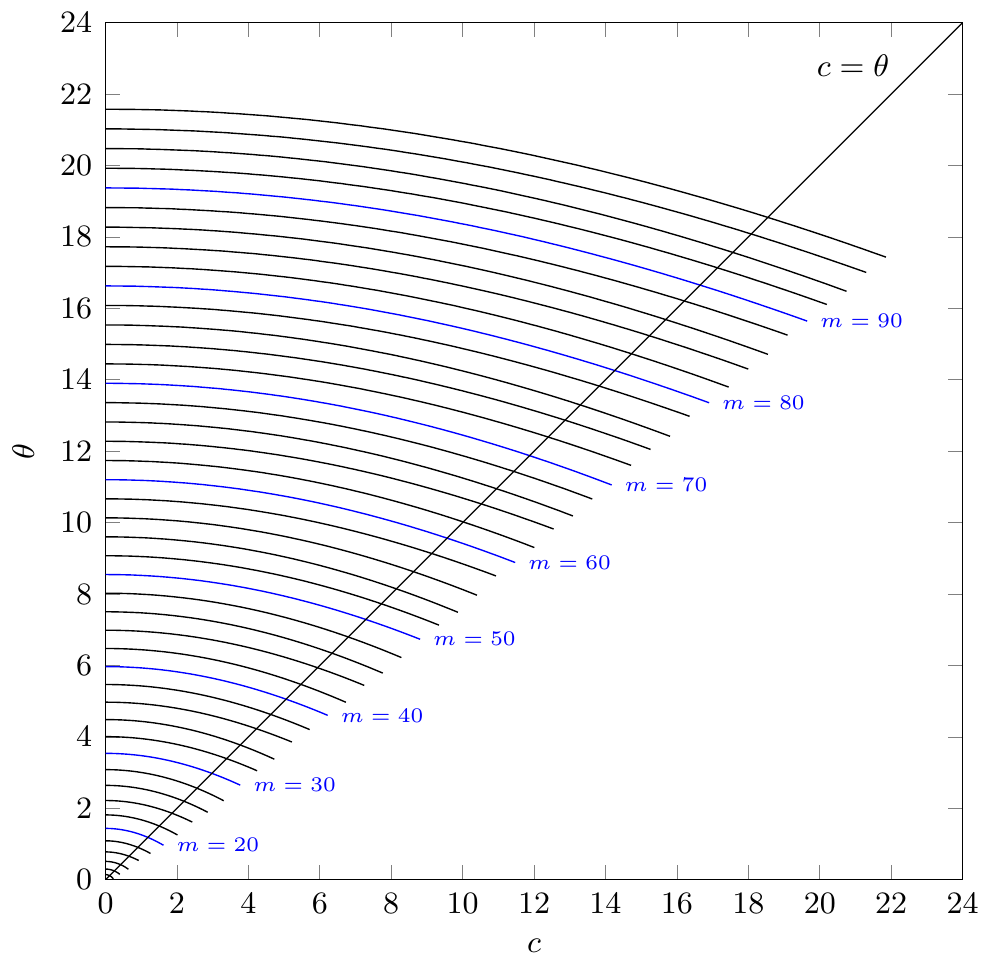}
	\caption{\label{fig:thetavsbeta_sym}The root $\theta=\theta_{m,\interval}$ as a function of the right endpoint $\interval$ of the interpolation interval for complex conjugate Leja points in $\rmi[-\interval,\interval]$. The tolerance is set to $\tol=2^{-53}$. Along each line the interpolation degree $m$ is kept fixed.}
\end{figure}

On the interval $\rmi\,[-\interval, \interval]$ on the imaginary axis, symmetrized or conjugate complex Leja points are defined as
\begin{align*}
	\xi_m\in\argmax_{\xi\in \rmi[-\interval,\interval]} \prod_{j=0}^{m-1}|\xi-\xi_j|, \quad \xi_{m+1}=-{\xi}_m\quad \text{for $m\geq1$ odd, and $\xi_0=0$}.
\end{align*}
We use conjugate complex pairs of points rather than standard Leja points in an interval along the imaginary axis as this allows real arithmetic for real arguments. 
This gives rise to polynomials of even degree. 

To allow conjugate complex Leja points in our backward error analysis we only need to change the actual computation of the values $\theta_{m,\interval}, \theta_m$ and $\gamma_{m,\interval}$. 
The theory itself stays the same. 
Figure~\ref{fig:thetavsbeta_sym} displays the path of $\theta_{m,c}$ for complex conjugate Leja points in $\rmi\,[-\interval,\interval]$. Table~\ref{tab:thetam:complex} displays a selection of (rounded) $\theta_m$ values for various tolerances.


If we apply complex conjugate Leja points in the framework of Section~\ref{sec:convergence_ellipses} we get ellipses for which the major axis is on the imaginary axis. 

\begin{table}\small\centering
\begin{tabular}{r|rrrrrrr}
 \hline
    $m$ &       10 &       20 &       30 &       40 &       50 \\
   half & 1.94e+00 & 4.53e+00 & 7.11e+00 & 9.62e+00 & 1.21e+01 \\
 single & 8.11e-01 & 2.99e+00 & 5.41e+00 & 7.85e+00 & 1.03e+01 \\
 double & 1.16e-01 & 1.19e+00 & 2.98e+00 & 5.06e+00 & 7.29e+00 \\
 \hline\hline
    $m$ &       60 &       70 &       80 &       90 &      100 \\
   half & 1.46e+01 & 1.70e+01 & 1.95e+01 & 2.20e+01 & 2.44e+01 \\
 single & 1.27e+01 & 1.52e+01 & 1.77e+01 & 2.01e+01 & 2.25e+01 \\
 double & 9.57e+00 & 1.19e+01 & 1.43e+01 & 1.67e+01 & 1.90e+01 \\\hline
\end{tabular}
\caption{\label{tab:thetam:complex}Samples of the (rounded) values $\theta_m$ with tolerances \emph{half} ($\tol=2^{-10}$), \emph{single} ($\tol=2^{-24}$) and \emph{double} ($\tol=2^{-53}$) for complex conjugate Leja interpolation.}
\end{table}

\subsection{\texorpdfstring{Extension to {$\varphi$} functions}{Extension to phi functions}}

The presented backward error analysis extends in a straightforward way to the $\varphi$ functions which play an important role in exponential integrators, see \citep{HO2010}.
We illustrate this here for the $\varphi_1$ function. For $A\in\mathbb{C}^{n\times n}$ and $w\in\mathbb{C}^n$ we observe that
\begin{align*}
	\varphi_1(A)w=[I,0]\exp\left( \begin{bmatrix}A & w\\0&0 \end{bmatrix}
	\right) \begin{bmatrix}0\\1\end{bmatrix},
\end{align*}
see \citep{Sidje1998}.
For the choice 
\begin{align*}
\mathcal{A}=\begin{bmatrix}A & w\\0&0 \end{bmatrix}\quad\text{and}\quad v=\begin{bmatrix}0\\1\end{bmatrix}
\end{align*}
the backward error is preserving the structure, i.e.
\begin{align*}
\Delta\mathcal{A}=\begin{bmatrix}\Delta A & \Delta w\\0&0 \end{bmatrix}.
\end{align*}
Thus the above analysis applies. 
For general $\varphi$ functions we can extend our approach with the help of \citep[Thm.~2.1]{Higham2011}.
\section{Additional aspects of interpolation}\label{sec:petc}

By using the previously described backward error analysis to compute the values $m_*$, $\nis_*$ and $c_*$ a working algorithm can be defined. 
Nevertheless, the performance of the algorithm can be significantly improved by some preprocessing and by introducing an early termination criterion. 
Moreover, interpolation in nonexact arithmetic will suffer from roundoff errors, in particular in combination with the hump phenomenon.
We address all these issues in this section.

\subsection{Spectral bounds and shift}\label{sec:shift}
In the above backward error analysis, it was assumed that the rectangle $R$ lies symmetrically about the origin.
In general, this requires a shift of $A$.
On the other hand, it is clear that a well chosen shift $\mu$ satisfying $\|A-\mu I\|\leq\|A\|$ is beneficial for the interpolation (a lower degree or less scaling steps will be required).
For the exponential function such a shift can easily be compensated by scaling. 
More precisely, if the shift $\mu$ is selected, we use 
\begin{align*}
	[\rme^{\mu/\nis}L_{m,\interval}(\nis^{-1}(A-\mu I))]^\nis
\end{align*}
as approximation of $\rme^{A}$. 

For our algorithm a straightforward shift is to center the rectangle $R$ at the origin, namely
\begin{align}\label{eq:shift}
	\mu=\frac{\SR+\LR}{2}+\rmi \frac{\SI+\LI}{2}.
\end{align}
If $A$ is real then $\SI=-\LI$ and $\mu\in\RR$. 
Therefore, a complex shift is only applied to complex matrices. 

It is easy to see that for a Hermitian or skew Hermitian matrix $A$ the proposed shift \eqref{eq:shift} coincides with the norm-minimizing shift presented in \citep[Thm.~4.21(b)]{Higham2008}. 
For a general matrix, the shift somewhat symmetrizes the spectrum of the matrix with regard to its estimated field of values. 

The shift $\mu=n^{-1}\operatorname{trace}A$ used in \citep{Higham2011} is a transformation that centers the spectrum of the matrix around the average eigenvalue. 
In many cases the two shifts are similar.
Nevertheless, it is possible to find examples where one shift leads to better results than the other. 
The matrix \texttt{one-sided} of Example~\ref{eg:longtimestep} is one of these cases.
For the trace shift a symmetrization of the rectangle $R$ might be required, resulting in a potential increase of scaling steps for the estimate based on \eqref{eq:ellipsestimate}.
For the method proposed here, we always use \eqref{eq:shift} as shift. 

\subsection{Early termination criterion}\label{sec:earlyterm}
The estimates based on \eqref{eq:thetam} and \eqref{eq:gammafoc} are worst case estimates and in particular do not take $v$ into account. 
As a result, the choice of $m_*$ is likely to be an overestimate and can be reduced in the actual computation. 
By limiting $m$ in the computation of $L_{m,\interval}(\nis^{-1}A)v^{(i)}$ in \eqref{eq:scalingsteps} we introduce a relative forward error that again should be bounded by the tolerance $\tol$. 
We propose to take
\begin{align}\label{eq:aposteriori}
	\|e_k\| &= \|L_{k,\interval}(\nis^{-1}A)v^{(i)} - L_{k-1,\interval}(\nis^{-1}A)v^{(i)}\| \nonumber\\
	&= \left|\exp[\xi_0,\ldots,\xi_k]\right|\left\| \prod_{j=0}^{k-1}(\nis^{-1}A-\xi_j I)v^{(i)}\right\|\leq \frac{\tol}{\nis} \|L_{k,\interval}(\nis^{-1}A)v^{(i)}\|
\end{align}
as an a posteriori error estimate for the Leja method in the $k$th step. Experiments show that \eqref{eq:aposteriori} turns out to be a good choice. 
In contrast to the method described in \citep{Higham2011} we divide the tolerance by the amount of scaling steps. 
This potentially increases the number of iterations per step but in practice results in a more stable computation for normal matrices, see Section~\ref{sec:numexp}.
Nevertheless, it sometimes leads to results of higher accuracy than prescribed.
In practice it is advisable to take the sum of two or three successive approximation steps for the estimate as this captures the behavior better.
On the other hand, it can also be beneficial to make the error estimate only every couple of iterations rather than in each step to save computational cost, see \citep{Leja-comp}. 
A second approach for an a posteriori error estimate for the Leja method based on the computation of a residual can be found in \citep{Leja-botchev}. 
This procedure can also be used here.
Furthermore, it is possible to adapt the early termination criterion to complex conjugate Leja points.
With the help of an early termination criterion computational cost can be saved for certain matrices, see Section~\ref{sec:numexp}.

\subsection{Handling the hump phenomenon}\label{sec:hump}

\begin{figure}[ht]\centering\captionsetup[subfigure]{aboveskip=-1pt}
	\begin{subfigure}{\linewidth}
	\includegraphics{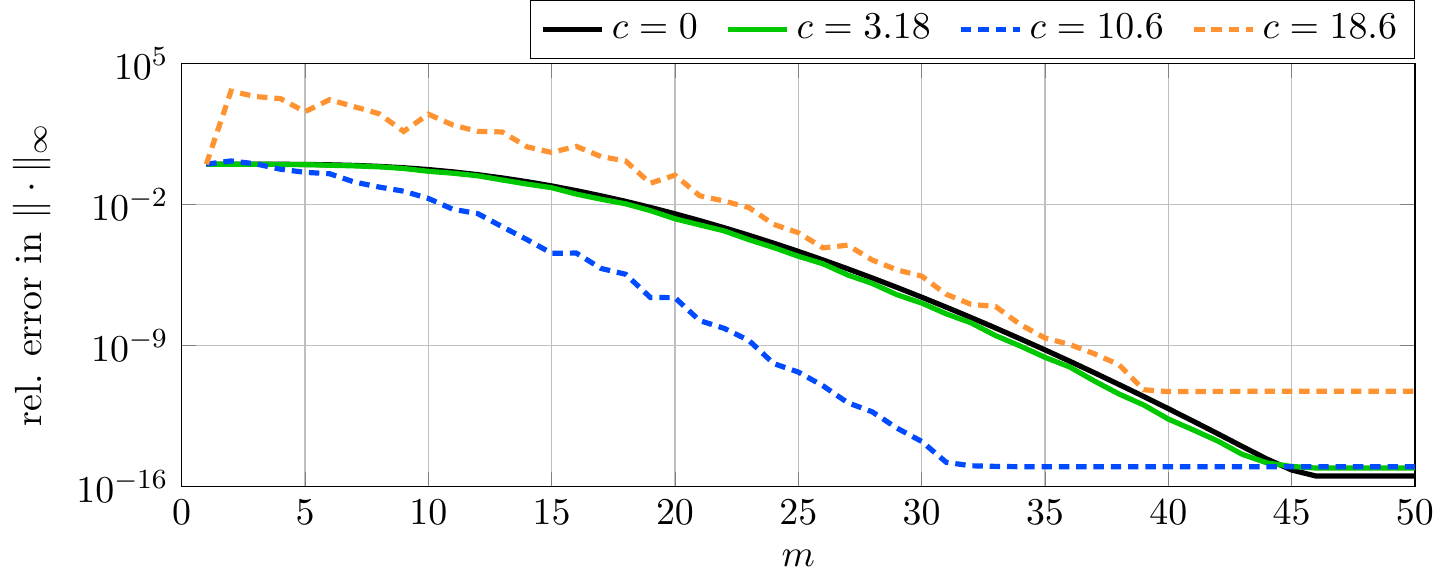}
	\caption{\label{fig:interror_a}Relative error vs.~degree of interpolation $m$ for $L_{m,\interval}(A)v$; real case.}
	\end{subfigure}\\[\abovecaptionskip]
	\begin{subfigure}{\linewidth}
	\includegraphics{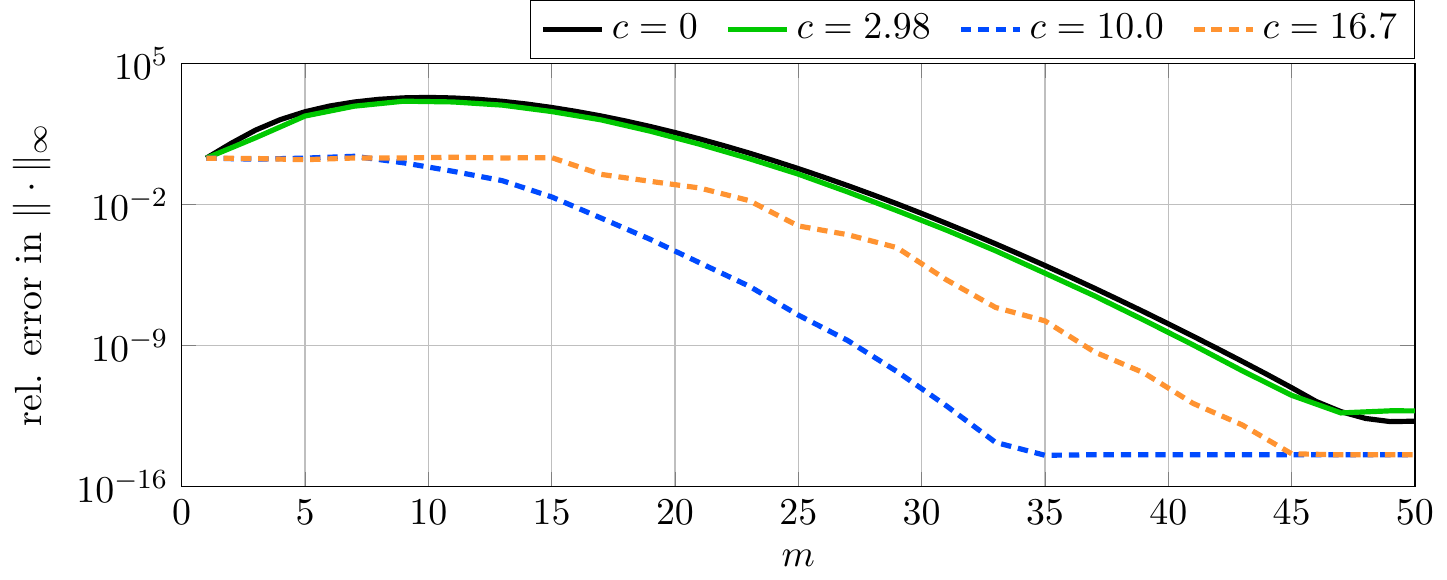}
	\caption{\label{fig:interror_b}Relative error vs.~degree of interpolation $m$ for $L_{m,\interval}(\mathrm{i}\,A)v$; complex case.}
	\end{subfigure}
	\caption{\label{fig:interror}Illustration of the hump phenomenon for the real and complex case. The used matrix $A=\text{\texttt{diag(linspace(-10,10,10))}}$ and $v=[1,\ldots,1]^\mathrm{T}$.}
\end{figure}

In general, the interpolation error does not decrease monotonically with the degree of interpolation. 
Even worse, a distinct hump can occur in certain situations, see Figure~\ref{fig:interror}.
This hump can significantly reduce the accuracy of the interpolation due to roundoff errors. 
The phenomenon is linked to the distribution of the eigenvalues of a matrix with respect to interpolation interval.
Note that the hump we are describing here is not the same as the one described in \citep{Moler2008} for nonnormal matrices. 

Figure~\ref{fig:interror} illustrates the problem for the matrix $A=\text{\texttt{diag(linspace(-10,10,10))}}$ and vector $v=[1,\ldots,1]^\mathrm{T}$. 
Figure~\ref{fig:interror_a} shows the real case. 
For $\interval=0$ (i.e.~the truncated Taylor series method) the final error is low and no hump is formed.
If we increase the interpolation interval $[-\interval,\interval]$, we observe that the necessary degree of the interpolation gradually decreases, while the final error stays approximately constant. 
The \emph{optimal} interpolation interval is reached when $\interval$ approaches the largest eigenvalue. 
In the figure this is the case for $\interval=10.6$.
When the interval is increased further, however, a hump starts to form.
This is due to the fact that the divided differences are significantly larger than the result, which has size $\rme^{10}$. 

As can be seen in Figure~\ref{fig:interror_b}, the behavior is different for the complex case. 
Here the divided differences have modulus one and a hump forms if the interpolation interval is too small. 
Note that the smallest necessary degree of interpolation is again obtained by selecting the \emph{optimal} interval.

In both cases the undesired behavior can be improved by obtaining a better estimate of the spectral radius and consequently reducing the interpolation interval. 
For this we employ the values $d_p=\|A^p\|^{1/p}$ which satisfy the well known relation 
\begin{align*}
\rho(A)=\lim_{p\to\infty}\|A^p\|^{1/p}.
\end{align*}
As long as the sequence of $\{d_p\}$ decreases, we adjust the interpolation interval accordingly. 

For a general matrix $A$ this phenomenon will influence the computation whenever $\|A\|$ strongly overestimates $\rho(A)$. 
In this case our algorithm chooses an interpolation interval that is far too large. Note that this happens in particular for nonnormal matrices. 

For the estimate based on \eqref{eq:thetam} the reduction of the interpolation interval is possible due to the behavior of the $\theta_{m,\interval}$ curve shown in Figure~\ref{fig:thetavsbeta}. 
However, if we use the estimate based on \eqref{eq:ellipsestimate} the reduction of the interpolation interval is not straightforward. 
If we fit our rectangle $R$ into an ellipse with semi-axis given by \eqref{eq:semiaxis} then, in general, $R$ will not fit into an ellipse with a smaller interpolation interval. 
We overcome this problem by adding a circle to the ellipses.
We use the largest circle defined by $a_{m,\theta_k}=b_{m,\theta_k}=\theta_k$ for some $k\leq m$ that fulfills \eqref{eq:gammafoc}, see Figure~\ref{fig:ellips} for an example. 
In most cases the radius of the circle is not going to be $\theta_m$.
If the values $d_p$ indicate a reduction of the interval, we restrict the ellipse estimate to these circles and perform a reduction of the interpolation interval. 
The validity of this process can be checked in the same manner as for~$\theta_m$.

\begin{rem}
	In the current version the algorithm does not allow to reduce the number $\nis$ along with the decay of $d_p$ as in \citep{Higham2011}. 
	Nevertheless, if a drastic decay is indicated it is possible to transform our method into Taylor interpolation by simply setting $\interval=0$.
\end{rem}

\section{Numerical examples} \label{sec:numexp}

In order to illustrate the behavior of our method we provide a variety of numerical examples. 
We use matrices resulting from the spatial discretization of time dependent partial differential equations already used in \citep{Leja-comp}. 
Furthermore, we also utilize examples from \citep{Higham2011} and certain prototypical examples to illustrate some specific behavior of the method. 
All our experiments are carried out with Matlab 2013a. 
As a measure of the required computational work we use the number of matrix-vector products (mv) performed by the method, without taking into account preprocessing tasks.
We will compare our method to the function \texttt{expmv} of \citep{Higham2011}. 

Note that the Leja method also employs divided differences. 
They are computed as described in \citep{Caliari2007}.
The used divided differences are precomputed as the employed interpolation intervals are fixed.

In the following we are going to compare different variations of our algorithm based on the presented ways to compute the scaling factor $\nis$ and degree of interpolation $m$. 
\begin{description}
	\item[\texttt{Algorithm 1}:] Uses $m_*$ and $\nis_*$ given by \eqref{eq:scalingstep_circle}. 
	\item[\texttt{Algorithm 2}:] Uses $m_*$ and $\nis_*$ given by \eqref{eq:scalingstep_ellipse}. 
\end{description}
In both algorithms, the early termination criterion \eqref{eq:aposteriori} is used, as well as the shift and the hump test discussed in the previous section, if not indicated otherwise. 
For \algre{} the hump test procedure is only employed if the estimate of the scaling step is based on circles. 

In Example~\ref{eg:lesp} we take a look at the stability of the methods with and without early termination, Examples~\ref{eg:ade} and \ref{eg:hump} focus on the selection of the degree and the interpolation interval for the different variations of our algorithm, and Example~\ref{eg:longtimestep} investigates the behavior for multiple scaling steps, i.e.,~$\nis>1$.

\begin{eg}[early termination]\label{eg:lesp} 
This example is taken from \citep[Exp.~2]{Higham2011} to show the influence of the early termination criterion for a specific problem.
We use the matrix $A$ as given by \texttt{gallery('lesp',10)}. 
This is a nonsymmetric, tridiagonal matrix with real, negative eigenvalues. 
We compute $\rme^{tA}v$ by \algnt{} and \algre{}, respectively, for 50 equally spaced time steps $t\in[0,100]$. We select the tolerance $\tol=2^{-53}$ and  $v_i=i$. 
As $A$ is a nonnormal matrix, \algre{} is restricted to circles to allow the hump reduction procedure.
\begin{figure}\centering
	\includegraphics{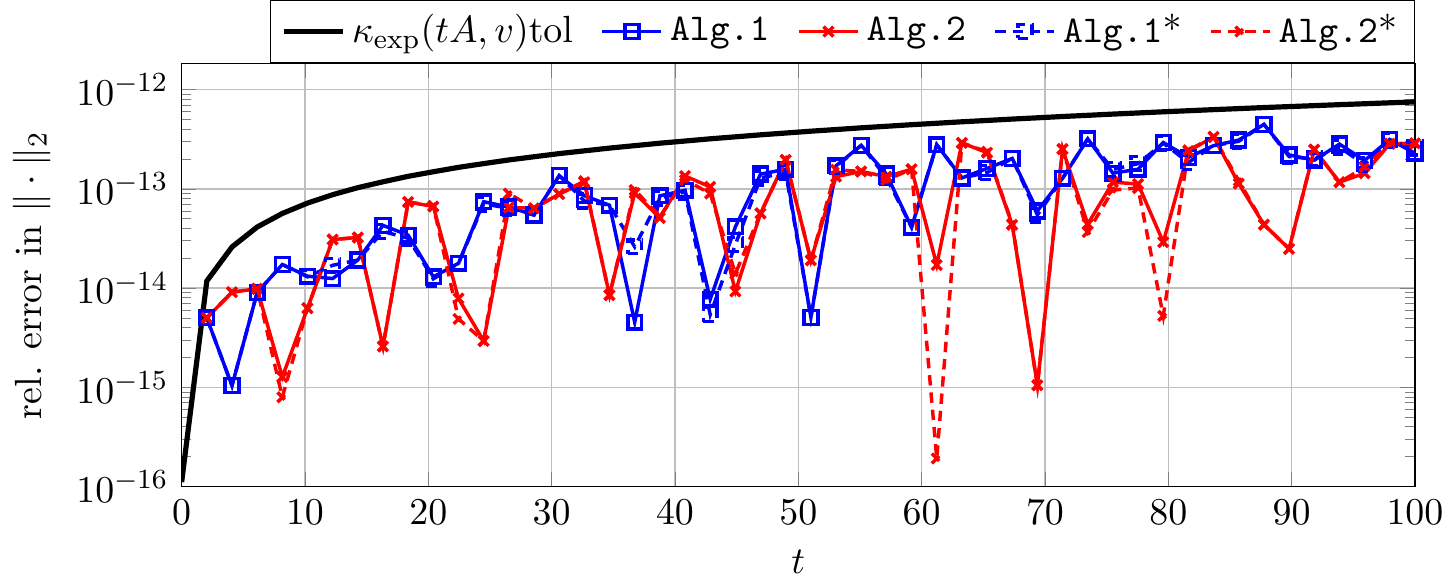}
	\caption{\label{fig:fs}Time step $t$ versus relative error in the $2$-norm for the computation of $\rme^{tA}v$ with tolerance $\tol=2^{-53}$. The $^*$ indicates that no early termination was used. Note that there is almost no visible difference between the methods with and without the early termination in place.}
\end{figure}
The results of the experiments can be seen in Figure~\ref{fig:fs} where the solid line corresponds to the condition number \eqref{eq:condnr} multiplied by the tolerance. 
We can not expect the algorithms to perform much better than indicated by this line. 
As condition number we use
\begin{align}\label{eq:condnr}
	\kappa_{\exp}(tA,v):=\frac{\|\exp(tA)\|_{2}\|v\|_{2}}{\|\exp(tA)v\|_{2}} + 
	\frac{\|(v^{\mathrm{T}}\otimes I)K_{\exp}(tA)\|_2\|\operatorname{vec}(tA)\|_{2}}{\|\exp(tA)v\|_{2}},
\end{align}
as defined in \citep[Eq.~(4.2)]{Higham2011}. Here $\operatorname{vec}$ denotes the vectorization operator that converts its matrix argument to a vector by traversing the matrix column-wise. Furthermore, let $L(A,\Delta A)$ denote the Fr\'{e}chet derivative of $\exp$ at $A$ in direction $\Delta A$ given by 
\begin{align*}
	\rme^{A+\Delta A}=\rme^A+L(A,\Delta A)+o(\|\Delta A\|).
\end{align*}
With the relation $\operatorname{vec}(L(A,\Delta A))=K_{\exp}(A) \operatorname{vec}(\Delta A)$ the Fr\'{e}chet derivative is given in its Kronecker form as $K_{\exp}(A)$. 
In addition we use the relation 
\begin{align*}
	\operatorname{vec} (L(A, \Delta A)v) = (v^{\mathrm{T}}\otimes I) \operatorname{vec} (L(A, \Delta A)).
\end{align*} 
For the computation we use the function \texttt{expm\_cond} from the Matrix Function Toolbox, see \citep{mft}. 

Overall we can see that the algorithms behave in a forward stable manner for this example. 
The early termination criterion shows no significant increase in the error. Both algorithms are well below $\kappa_{\exp}(tA,v)$ for all values of $t$.
For this rather small matrix we used the exact norm and not a norm estimate to allow for a better comparison.
\end{eg}

\begin{eg}[advection-diffusion equation]\label{eg:ade}
In order to show the difference between the two suggested processes for selecting the interpolation interval for our algorithms, we use an example that allows us to easily vary the spectral properties of the discretization matrix. 
Let us consider the advection-diffusion equation 
\begin{align*}
	\partial_t u=a \Delta u + b (\partial_x u+\partial_y u) 
\end{align*} on the domain $\Omega=[0,1]^2$ with homogeneous Dirichlet boundary conditions. 
This problem is discretized in space by finite differences with grid size $\Delta x =(N+1)^{-1}$, $N\ge1$. 
As a result of the discretization we get a sparse $N^2\times N^2$ matrix $A$. 
We define the grid P\'{e}clet number
\begin{align*}
	\Pe=\frac{|b| \Delta x}{2a}
\end{align*} as the ratio of advection to diffusion, scaled by $\Delta x$. 
By increasing $\Pe$ the nonnormality of the discretization matrix can be controlled. 
In addition, $\Pe$ describes the height-to-width ratio of the rectangle $R$. 
The estimates for $\SR$ and $\LR$ stay constant but $\SI=-\LI$ increases with $\Pe$. 

For the following computations the parameters are chosen as: $N=20$, $a=1$ and $b=\tfrac{2a\Pe}{\Delta x}$. 
As a result, for $\Pe=0$ we get that $R$ is an interval on the real axis and for $\Pe=1$ a square.
For $\Pe=0$ the matrix is equal to $\texttt{-(N+1)\^{}2*gallery('poisson',N)}$. 
The vector $v$ is given by the discretization of the initial value $u_0(x,y)=256\cdot x^2 (1-x)^2 y^2 (1-y)^2$.
In the following discussion we call the shifted matrix again $A$.

\begin{table}\small\centering
\noindent\begin{tabular}{l|*{10}{c|}c}
 			& \multicolumn{4}{|c|}{\algnt}   				& \multicolumn{4}{|c|}{\algre}				& \multicolumn{3}{|c}{\expmv}	\\
 			& $m_*$	&$m$	&	rel.~err 	&$\interval$	&$m_*$	&$m$	&	rel.~err &$\interval$	& $m_*$	& $m$	&	rel.~err 	\\\hline
$\Pe=0$		&	54	&	32	&	3.66e-15	&	8.96		& 40 & 32 & 3.66e-15 &  8.96			&	52	&	44	&	2.88e-15	\\
$\Pe=0.2$	&	54	&	34	&	5.47e-15	&	8.96		& 49 & 35 & 4.17e-15 &  8.28			&	52	&	44	&	3.91e-15	\\
$\Pe=0.4$	&	54	&	35	&	2.21e-15	&	8.96		& 55 & 36 & 2.86e-15 &  9.24			&	52	&	44	&	2.34e-15	\\
$\Pe=0.6$	&	54	&	38	&	3.63e-15	&	8.96		& 62 & 40 & 6.36e-15 & 11.08			&	52	&	44	&	4.93e-15	\\
$\Pe=0.8$	&	54	&	41	&	2.98e-15	&	8.96		& 67 & 43 & 9.86e-15 & 11.34			&	52	&	43	&	2.59e-15	\\
$\Pe=1$		&	54	&	44	&	3.24e-15	&	8.96		& 72 & 49 & 4.50e-14 & 12.99			&	52	&	39	&	1.30e-15		
\end{tabular}
\caption{\label{tab:ad_singlestep}For varying grid P\'{e}clet number in Example~\ref{eg:ade} the selection of the degree of interpolation $m_*$, the actual degree due to the early termination $m$ and the right endpoint $\interval$ of the interpolation interval are shown. We compute $\exp(tA)v$ with a time step $t=5$e-3, discretization parameter $N=20$ and tolerance $\tol=2^{-53}$.  The error is measured relative to the result of the Matlab built-in function \texttt{expm} in the maximum norm.}
\end{table}

Table~\ref{tab:ad_singlestep} gives the results of an experiment where we varied the grid P\'{e}clet number. 
We show the results of the different selection procedures. 
The time step $t=5$e-3 is chosen such that \expmv{} is able to compute the result without scaling the matrix. 
The actual degree of interpolation $m$ and the relative error with respect to the method \texttt{expm} in the maximum norm are shown. 
As the maximum norm of the matrix stays the same (for fixed $N$) the parameters of \expmv{} and \algnt{} are always the same. 
For $\Pe=1$ the eigenvalues of the matrix $tA$ are in a small circle around zero and therefore the Taylor approximation requires a lower degree of interpolation.

On the other hand we can see that for a small height-to-width ratio (small $\Pe$) the estimate based on ellipses, i.e.~\algre{}, produces a significantly smaller $m_*$ with the same actual degree $m$ of interpolation and comparable error. 
This means that less scaling is required for larger $t$, cf.~Example~\ref{eg:longtimestep}.
When the rectangle $R$ is closer to a square the algorithm still produces reliable results but is slightly less efficient than \algnt{}. 

\begin{rem}[$\theta_m$ selection]\label{rem:thetamax}
As mentioned in Section~\ref{sec:cd} it is possible to select the $\theta_m$ values differently depending on $\interval$. By selecting $\hat\theta_m=\max_\interval\theta_{m,\interval}$ the computation corresponding to $\Pe=0$ gives the following results: $m_*=51$, $m=44$ and $\interval=4.31$. This indicates a slower convergence with a similar error, as the eigenvalues of the shifted matrix are in $[-8.82,8.82]$ but we interpolate in $[-4.31,4.31]$.
\end{rem}
\end{eg}

\begin{eg}[hump and scaling steps] \label{eg:hump}
In order to illustrate the potential gain of testing for a hump in our algorithm, we use the matrix $A$ given by \texttt{-1/2*gallery('triw',20,4)} and $v_i=\cos\,i$.
This corresponds to \citep[Experiment 6]{Higham2011} with a single time step of size $t=0.5$. 
In the following discussion we call the shifted matrix again $A$.

The $20\times20$ matrix $A$ is an upper triangular matrix with $0$ in the main diagonal and $-2$ on the strict upper triangular part.
The $1$-norm of the matrix is $\|A\|_1= 38$ and $\rho(A)=0$.  
For this example the truncated Taylor series method is optimal and stagnates after $20$ iterations, in a single scaling step, with a final error of about $10^{-14}$. 
We use this example to illustrate the hump phenomenon and the procedure to counteract it.
This will result in a better performance of the Leja method, even though for this example it is not as efficient as \expmv{}.
\begin{figure}\centering
	\includegraphics{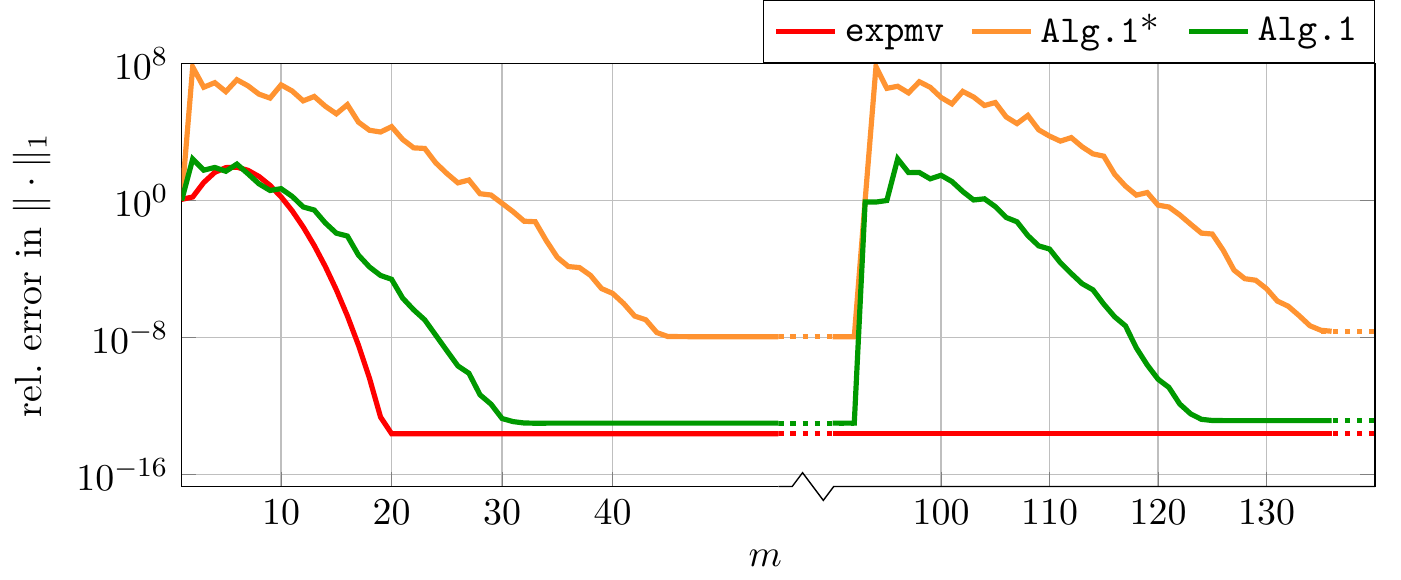}
	\caption{\label{fig:triw_hump}Relative error vs.~interpolation degree $m$ for the approximation of $\exp(A)v$. The plot illustrates the behavior of the hump reduction procedure. Here \algnt{}* indicates that no hump reduction was performed. The relative error is measured with respect to $\exp(A/2)v$ in the first scaling step and $\exp(A)v$ in the second. No early termination criterion is used in the computation.}
\end{figure}

In Figure~\ref{fig:triw_hump} we illustrate the behavior of the hump reduction procedure, described in Section~\ref{sec:petc}. 
A first observation is that our method selects two scaling steps ($\nis=2$). 
For this nonnormal matrix the rectangle $R$ is a square and therefore we only show the results for \algnt{}, as we can expect a better performance for this algorithm (cf.~Example~\ref{eg:ade}).
We can see quite clearly that a hump of about $8$ digits is formed when we use the initial guess of $\interval=\theta_{92}=19.10$. 
As expected we get an error of about $10^{-8}$ in each scaling step and consequently a final error of the same size. 
In this experiment we deactivate the early termination criterion and therefore the algorithm uses $m=92$ in each of the two scaling steps. 
Furthermore, for \algnt{}* and \algnt{} the relative error is measured with respect to $\exp(A/2)v$ in the first scaling step and $\exp(A)v$ in the second.
For sake of comparison we run \expmv{} without early termination and $182$ iterations, and we measure the relative error with respect to $\exp(A)v$.

On the other hand, if we reduce the interval length, the hump is reduced as well. 
For this (shifted) matrix the values $d_p$ are 
$(38, 26, 20, 16, 13, \ldots, 0)$ 
(see Section~\ref{sec:hump}), where $d_{20}=0$. 
These values suggest that the interpolation interval should be reduced. 
In fact, by reducing the interval to $\interval=0$ we would recover the truncated Taylor series method and therefore the optimal choice for this example.
The cost of the computation of $d_p$ can not be neglected and in a practical implementation $p$ is therefore limited.
Experiments have shown that $p\leq 5$ is a practical choice. 
Therefore, we use the interpolation interval $\interval=\theta_{45}=6.67$ in the reduced case. 

From Figure~\ref{fig:triw_hump} we can see that the algorithm to reduce the impact of the hump is working. 
The hump is significantly reduced and the error is close to the error of the truncated Taylor series method.
Nevertheless, the method takes about three times as many iterations (approximately $60$) than the truncated Taylor series method, cf.~Example~\ref{eg:longtimestep}. 

Even though this procedure might not be necessary for accuracy it is still beneficial for the overall cost reduction. 
If we only require $8$ digits of accuracy we do not need to reduce the interpolation interval to achieve this. 
However, for this example it would still be beneficial to reduce the interpolation interval, as the necessary degree of interpolation is reduced as well. 
This is related to the observations of Remark~\ref{rem:thetamax} and Section~\ref{sec:hump}. 
Here the norm of the matrix is a large overestimate of the spectral radius leading to slow convergence. 
\end{eg}


\begin{eg}[behavior for multiple scaling steps]\label{eg:longtimestep} 

In this experiment we investigate the behavior of the methods for multiple scaling steps.
We use the two matrices of Examples~\ref{eg:ade} and~\ref{eg:hump} from above and in addition the matrices \texttt{orani676} and \texttt{bcspwr10} which are obtained from the University of Florida sparse matrix collection \citep{Davis2011}, as well as several other matrices, see Table~\ref{tab:matrixprob}.
The sparse matrix \texttt{orani676} is real and nonnormal with $90158$ nonzero entries,
whereas \texttt{bcspwr10} is a real and symmetric sparse matrix with $13571$ nonzero entries.
The matrix \texttt{triu} is an upper triangular matrix with entries uniformly distributed on $[-0.5,0.5]$.
For the matrix \texttt{onesided} we have a $41\times41$ upper triangular matrix with one eigenvalue at $10$ and $40$ eigenvalues uniformly distributed on $[-10.1,-9.9]$, with standard deviation of $0.1$.
The values in the strict upper triangular part are uniformly distributed on $[-0.5,0.5]$.
Furthermore, we use \texttt{S3D} from \citep[Example 3]{Leja-comp}, a finite difference discretization of the three dimensional Schr\"odinger equation with harmonic potential in $[0,1]^3$. The matrix \texttt{Trans1D} is a periodic, symmetric finite difference discretization of the transport equation in $[0,1]$.
For the matrices \texttt{orani676}, \texttt{S3D} and \texttt{Trans1d} complex conjugate Leja points are used in the computation.

As vector $v$ we use $[1,\ldots,1]^\mathrm{T}$ for \texttt{orani676}, $[1,0,\ldots,0,1]^\mathrm{T}$ for \texttt{bcspwr10}, $v$ as specified in Example~\ref{eg:ade} for \texttt{AD}, the discretization of $4096x^2(1-x)^2y^2(1-y)^2z^2(1-z)^2$ is used for \texttt{S3D}, the discretization of $\exp(-100(x-0.5)^2)$ for \texttt{Trans1D}, and $v_i=\cos\,i$ for all other examples. 
This corresponds to \citep[Exp.~7]{Higham2011}. 
\begin{table}[ht]\small
\begin{subtable}{\linewidth}\centering\setlength{\tabcolsep}{4pt}
\begin{tabular}{l|l|*{4}{c}l*{3}{c}}
\# &\multicolumn{1}{c|}{$A$}			&	$n$	& $t$		& $\LR-\SR$		& $\LI-\SI$ & \multicolumn{1}{c}{$\kappa_1$} & $\|\cdot\|_1$ & $\|\cdot\|_2$ & $\|\cdot\|_\infty$ \\\hline
1 & orani676 & 2529 &  100  & 1.0e+03 &  1.0e+03 & 0.002   & 1.0e+03 & 3.2e+01 & 9.4e+00 \\
2 & bcspwr10 & 5300 &   10  & 2.6e+01 &  0.0e+00 & 0       & 1.4e+01 & 6.8e+00 & 1.4e+01 \\
3 &     triw & 2000 &   10  & 8.0e+03 &  8.0e+03 & 0.5     & 8.0e+03 & 5.0e+03 & 8.0e+03 \\
4 &     triu & 2000 &   40  & 1.0e+02 &  1.0e+03 & 0.021   & 1.0e+03 & 4.2e+01 & 1.0e+03 \\
5 &  AD Pe=0 & 9801 &  1/4  & 8.0e+04 &  0.0e+00 & 0       & 8.0e+04 & 7.9e+04 & 8.0e+04 \\
6 &  AD Pe=0 & 9801 &    1  & 8.0e+04 &  0.0e+00 & 0       & 8.0e+04 & 7.9e+04 & 8.0e+04 \\
7 & onesided &   41 &    5  & 3.1e+01 &  1.2e+01 & 0.5     & 2.0e+01 & 1.1e+01 & 2.0e+01 \\
8 & S3D      &27000 &  1/2  & 0.0e+00 &  5.7e+03 & 0       & 5.8e+03 & 5.7e+03 & 5.8e+03 \\
9 & Trans1D  & 1000 &    2  & 0.0e+00 &  2.0e+03 & 0       & 1.0e+03 & 1.0e+03 & 1.0e+03 \\
\end{tabular}
\caption{\label{tab:matrixprob}Summary of the spectral properties of the matrices.}
\end{subtable}\\[5pt]

\begin{subtable}{\linewidth}\centering\setlength{\tabcolsep}{4pt}
\noindent\begin{tabular}{lr|*{9}{c|}c}
\multicolumn{2}{c|}{} & \multicolumn{3}{|c|}{\algnt{} $1$-norm}   	& \multicolumn{3}{|c|}{\algre{} $1$-norm}	& \multicolumn{3}{|c}{\expmv{} $1$-norm}			 \\
 \#	&	$t$& $\nis$&   $\mv$ & rel.err  & $\nis$ &     $\mv$ & rel.err  & $\nis$ &     $\mv$ & rel.err  \\\hline
1 &   100 & 4639 & 41751 & 1.8e-11 & 3508 & 31572 & 2.2e-11 &   21 &    526 & 4.0e-08 \\
2 &    10 &    6 &   157 & 7.8e-10 &    6 &   157 & 7.8e-10 &    5 &    171 & 7.2e-07 \\
3 &    10 & 3408 & 98840 & 5.4e-09 & 2423 & 87233 & 1.0e-07 & 1588 &  10425 & 1.1e-09 \\
4 &    40 & 1730 & 22495 & 1.6e-11 & 1268 & 17755 & 1.6e-12 &   59 &    960 & 4.3e-09 \\
5 &   1/4 &  427 & 14945 & 1.0e-08 &  357 & 13923 & 1.9e-09 &  749 &  29211 & 2.2e-06 \\
6 &     1 & 1705 & 59675 & 1.9e-08 & 1426 & 55614 & 3.3e-09 & 2995 & 116805 & 9.0e-06 \\
7 &     5 &    5 &   129 & 3.0e-09 &    4 &   119 & 1.6e-10 &    8 &    296 & 2.1e-08 \\
8 &   1/2 &   65 &  3185 & 2.2e-11 &   57 &  2793 & 8.0e-09 &  108 &   5400 & 1.3e-07 \\
9 &     2 &   89 &  4539 & 9.5e-13 &   79 &  3871 & 1.4e-08 &  150 &   5135 & 3.6e-08 \\
\end{tabular}
\caption{\label{tab:longtimestep}Results for each matrix and the used algorithms, respectively.}
\end{subtable}
\caption{\label{tab:longtime:everything}Results for Example~\ref{eg:longtimestep}. For a tolerance of $\tol=2^{-24}$ we compute $\exp(t A)v$ in a single call of the respective algorithm. The value $\nis$ indicates the number of scaling steps and $\mv$ denotes the number of matrix-vector products without preprocessing. The values $\SR$, $\LR$, $\SI$, $\LR$ correspond to \eqref{eq:recvalues}.}
\end{table}

We summarize the properties of all the matrices used in this example in Table~\ref{tab:matrixprob}. The tolerance is chosen as $2^{-24}$ and the relative error is computed with respect to \expmv{} running with the highest accuracy.
Furthermore we use 
\begin{align*}
\kappa_1=\frac{\|AA^*-A^*A\|_1}{\|A\|_1^2}
\end{align*}
as an indicator for the nonnormality of the matrices. From now on we refer to the matrices by their number given in the first column of Table~\ref{tab:matrixprob}.

We can see that for the nonnormal matrices $\{1,3,4\}$ the algorithm \expmv{} is superior in terms of matrix-vector products, in comparison to both variants of our algorithm. 
This is largely due to the fact that for these matrices the method \expmv{} can reduce the number of scaling steps based on the values $d_p$ (see Remark~\ref{rem:alphap}).
As the Leja method is not able to do this, the only way of getting comparable results for these example is by obtaining sharper bounds for the rectangle $R$ in \algre{}.
This could be achieved using the Matlab routines eig (based on LAPACK --- Linear Algebra PACKage and suited for full matrices), eigs (based on ARPACK --- ARnoldi PACKage and suited for sparse matrices), or an eigensolver of your choice fitted to the example.
However, the computation can be very expensive and therefore is not practicable in a general purpose algorithm. 

Furthermore, for these matrices the user specified norm has a relevant influence on the performance, as can be seen in Tables~\ref{tab:longtimestep} and~\ref{tab:longtimestep2}, respectively.
If the problem is considered with the $2$- or the maximum norm the number of matrix-vector products is significantly reduced. 

\begin{table}[t]\small\centering\setlength{\tabcolsep}{5pt}
\noindent\begin{tabular}{lr|*{6}{c|}c}
\multicolumn{2}{c|}{} & \multicolumn{3}{|c|}{\algnt{} $2$-norm} & \multicolumn{3}{|c|}{\algnt{} $\infty$-norm}	\\
 {\#}	&	$t$& $\nis$&   $\mv$ & rel.err   & $\nis$&   $\mv$ & rel.err  \\\hline
1 &   100 &  142 &  2434 & 1.8e-10 &    43 &  2195 & 2.2e-11 \\
2 &    10 &    2 &   106 & 7.0e-08 &     6 &   154 & 1.0e-09 \\
3 &    10 & 2175 & 76141 & 5.8e-08 &  3408 & 98847 & 5.4e-09 \\
4 &    40 &   66 &  2122 & 7.4e-10 &  1748 & 22732 & 1.5e-11 \\
7 &     5 &    3 &    91 & 5.5e-10 &     5 &   128 & 3.0e-09 \\
\end{tabular}
\caption{\label{tab:longtimestep2}For a tolerance of $\tol=2^{-24}$ we compute $\exp(t A)v$ in a single call of the algorithm \algnt{} with the $2$- and maximum norm, respectively. The value $\nis$ indicates the number of scaling steps and $\mv$ denotes the number of matrix-vector products without preprocessing. The numbers \# correspond to Table~\ref{tab:matrixprob}.}
\end{table}

For the matrices $\{2, 5, 6, 7, 8, 9\}$ the results show a different picture. 
Here, the Leja method is clearly beneficial in terms of matrix-vector products. 
Furthermore, we also produce a smaller error in comparison to the truncated Taylor series approach.
This is due to the fact that we divide the tolerance by $\nis$ in the early termination criterion, cf.~\eqref{eq:aposteriori}. 
In the case of the complex conjugate Leja points this leads to a higher accuracy than required.
On the other hand, for the \texttt{AD} problem, the errors of the Leja methods increase by a factor of $2$, if we increase $t$ by a factor of $4$, whereas the error for \expmv{} increases by a factor of $4$.
This is due to the fact that the used early termination criterion \eqref{eq:aposteriori} for the Leja method takes the number of scaling steps into account whereas \expmv{} does not. 

For matrices $\{1,8,9\}$ conjugate complex Leja points are used, see Section~\ref{sec:complexleja}. 
For the two normal matrices $\{8,9\}$ the Leja method saves a lot of matrix-vector products in comparison to \expmv{}. 
As here the rectangle $R$ is a line, \algre{} is again superior to \algnt{} as it leads to fewer scaling steps and less matrix-vector products. 

Matrix $2$ is normal. 
However, only for the $2$-norm we have that $\|A\|=\rho(A)$. 
This is the reason why the number of scaling steps is only two in the $2$-norm.

The final error of the methods is always comparable.
The more precautious approach we propose leads sometimes to an increase in accuracy.   
Nevertheless, in the cases where the Leja method is beneficial it still uses significantly less matrix-vector products than \expmv{}.

A comparison of \algnt{} and \algre{} shows that none of the two approaches can be considered superior or the \emph{better} overall choice. 
Due to the construction, \algre{} provides a scaling factor and a degree of interpolation that are independent of the norm, even though the Gerschgorin discs are closely related to the $1$- and maximum norm.
Nevertheless, the reduction of the interpolation interval is connected to a norm.
In fact, this is also the case where the two methods do provide similar estimates for $\nis$.
In total, if we always select the method with the least (predicted) computational cost we always use the more efficient methods as we save matrix-vector products.
This indicates that a combination of the two algorithms, where we always select the one with the least expected cost is beneficial. 

Depending on the specified norm, \algnt{} has some significant fluctuations in performance. 
\end{eg}

\section{Discussion}\label{sec:discussion}
The backward error analysis presented in this work provides a sound basis for the selection of the scaling parameter $\nis_*$ and the degree of interpolation $m_*$ for the Leja method. 
With this information at hand the algorithm becomes in a sense direct, as the maximal number of matrix-vector products is known after the initial preprocessing.
The cost of \algnt{} is determined by the norm of the matrix, whereas the cost of \algre{} is determined by the spectral information of the matrix. 
The convergence is monitored by the early termination criterion. 
The practical use of this approach is confirmed by the numerical experiments of Section~\ref{sec:numexp}.

The algorithm can be adapted in a similar way as the \expmv{} method to support dense output and provides essentially the same properties as \citep[Algorithm 5.2]{Higham2011}. In particular this means that the new algorithm also has some benefits in compared to Krylov subspace methods. 

Note that in certain applications one has to compute $\rme^{tA}V$ for a scalar $t$ and a $n\times n_0$ matrix $V$.  
This problem, however, is not more general since the product $tA$  can always be considered as a new matrix and the performed analysis extends to a matrix $V$ instead of a vector $v$. 
This is especially interesting in comparison to Krylov subspace methods as the available implementations would need to be called repeatedly for each column of $V$. 

In comparison to \expmv{} the Leja method is especially beneficial for matrices where the values $d_1,\ldots, d_p$ do not vary much. 
In these cases the method saves matrix-vector products. 
On the other hand our method makes a higher preprocessing effort than \expmv{}. 
This is due to the more complex selection procedure of $m_*$ and $\nis_*$ and the fact that we need an estimate of the field of values. 
As the overall cost is dominated by the matrix-vector products, this fact comes only into play for low-dimensional examples.

The combination of the two algorithms \algnt{} and \algre{}, where we select the scaling parameter and the degree of interpolation based on the minimum of the predicted cost of the two algorithms, seems to be the logical choice for a combined (black box) algorithm. With the changes applied to the method it can be called for any matrix $A$, it is numerically stable, the costs are predictable and the effort for the implementation is manageable.

In the present version of our algorithm, the knowledge of $d_p$ can not be used to properly scale the interpolation interval. However, and this is the focus of our future work, it is possible to modify the method and repeatedly use zero as interpolation point. 
These so-called Leja--Hermite methods will then be able to make use of $d_p$ in a suitable fashion as \expmv{}.

A Matlab implementation of the algorithm presented in this paper is available on the homepage \url{https://numerical-analysis.uibk.ac.at/exponential-integrators}.

\section*{Acknowledgment}
We thank the referees for their constructive remarks which helped us to improve the presentation of the paper considerably.
\bibliographystyle{plainnat}

\end{document}